\documentclass{article}
\usepackage{indentfirst,amssymb,amsmath}
\numberwithin{equation}{section}
\newtheorem{definition}{Definition}
\newtheorem{thm}{Theorem}[section]
\newtheorem{lem}{Lemma}
\newtheorem{prop}{Proposition}
\numberwithin{definition}{section}
\numberwithin{thm}{section}
\numberwithin{lem}{section}
\numberwithin{prop}{section}
\usepackage[left=2.0cm,right=2.0cm,top=2.5cm,bottom=1.5cm]{geometry}
\usepackage{hyperref}
\begin{document}
	\begin{center}
		{\Large{\bf Regularity results to the class of variational obstacle problems with variable exponent}}\vskip 4mm
		\large Debraj Kar\footnote{Department of Mathematics, University of Kalyani, India , Email : kardebraj98@gmail.com}
	\end{center}\vskip 2mm
	\begin{abstract}
		In this paper we are going to prove local gradient estimates and higher differentiability result for the solutions of variational inequalities 
		\begin{equation*}
			\int_\Omega\big<\mathcal{A}(x,u,Du),D(\phi-u)\big>dx\geq \int_\Omega\mathcal{B}(x,u,Du)(\phi-u)dx
		\end{equation*}
		for all $\phi\in \mathcal{K}_\psi(\Omega)$. Here $\Omega(\subset\mathbb{R}^n)$ is bounded , $n\geq 2$ and $\psi:\Omega\rightarrow\mathbb{R}$ is called obstacle. Here we deal with variable exponent growth , namely $p(.)$-growth . At first we prove Calder\'on-Zygmund estimate and then using this result to prove higher differentiability result in Besov scale. 
	\end{abstract}\vskip 1mm
	{\bf Keywords :} Obstacle problems, Variational inequalities, Variable exponent, Calder\'on-Zygmund estimates, Higher Differentiability.
	
	\section{Introduction}
	The aim of this paper to study Calder\'on-Zygmund estimate and Higher Differentiability in Besov scale of the solutions $u\in \mathcal{K}_\psi(\Omega)$ of variational obstacle problem of the following form
	\begin{equation}\label{eq1.1}
		\min\Bigg\{\int_\Omega f(x,v,Dv	) : v\in \mathcal{K}_\psi(\Omega)\Bigg\},
	\end{equation}
	where $\Omega$ is bouded subset of $\mathbb{R}^n$, for $n\geq 2$ and $f(x,u,\xi):\Omega\times \mathbb{R}\times\mathbb{R}^n\rightarrow\mathbb{R}$ is a strickly convex Carath\'eodory function with $\mathcal{C}^2$ class function with respect to $\xi$.\vskip 1mm
	The class $\mathcal{K}_\psi(\Omega)$ of admissible functions defined as follows 
	\begin{equation}\label{eq1.2}
		\mathcal{K}_\psi(\Omega)=\Big\{v\in W^{1,p(.)}(\Omega):v\geq \psi\;a.e.\;in\;\Omega\Big\}
	\end{equation}
	where $\psi:\Omega\rightarrow [-\infty,+\infty)$ is fixed function, called obstacle belongs to $W^{1,p(.)}(\Omega)$. Without any ambiguity I take that $\mathcal{K}_\psi(\Omega)$ is non-empty as $\psi$ itself in $\mathcal{K}_\psi(\Omega)$.\vskip 1mm
	Regularity work of solution of such problem have been done earlier see \cite{MZ49_43p_26,G20_15_1_24,EN43_14_1_14,EH50_6_10} etc .
	$u\in W^{1,p(x)(\Omega)}$ will be solution in $\mathcal{K}_\psi(\Omega)$ of the obstacle problem (\ref{eq1.1}) iff $u\in\mathcal{K}_\psi(\Omega)$ satisfies the following inequality
	\begin{equation}\label{eq1.3}
		\int_\Omega\big<\mathcal{A}(x,u,Du),D(\phi-u)\big>dx\geq \int_\Omega\mathcal{B}(x,u,Du)(\phi-u)dx
	\end{equation}
	where for all $\phi\in \mathcal{K}_\psi(\Omega)$ and  $\mathcal{A}(x,y,z)=D_zF(x,y,z)$, $\mathcal{B}(x,y,z)=-D_yF(x,y,z)$.\vskip 1mm
	For sake of simplicity we will take assumptions on $\mathcal{A}$ instead of $f$. we assume there exists $\nu,L,l$ and $p(x)\geq 2$ following conditions holds
	\begin{equation}\tag{A1}\label{eqa1}
		\big<\mathcal{A}(x,u,\xi)-\mathcal{A}(x,u,\eta),\xi-\eta\big>\geq \nu|\xi-\eta|^2(1+|\xi|^2+|\eta|^2)^\frac{p(x)-2}{2}
	\end{equation}
	\begin{equation}\tag{A2}\label{eqa2}
		|\mathcal{A}(x,u,\xi)-\mathcal{A}(x,u,\eta)|\leq L|\xi-\eta|(1+|\xi|^2+|\eta|^2)^\frac{p(x)-2}{2}
	\end{equation}
	\begin{equation}\tag{A3}\label{eqa3}
		|\mathcal{A}(x,u,\xi)|\leq l(1+|\xi|^2)^\frac{p(x)-1}{2}
	\end{equation}
	For almost every $x\in\Omega$ and every $u\in \mathbb{R}$, $\xi,\eta\in\mathbb{R}^n$. Moreover we suppose that there exist a  measurable function $a:\Omega\rightarrow [0,+\infty)$ and a positive exponent $r(x)<p(x)-1$ (see for reference \cite{G20_15_1_24}) such that
	\begin{equation}\tag{A4}\label{eqa4}
		|\mathcal{B}(x,\eta,\xi)|\leq |\xi|^{r(x)}+|\eta|^{r(x)}+a(x)
	\end{equation}
	For almost every $x\in\Omega$ and every $u\in \mathbb{R}$, $\xi,\eta\in\mathbb{R}^n$.\vskip 1mm
	Moreover we would like to assume an  assumption given following which will be important in sequel 
	\begin{equation}\tag{H1}\label{eqh1}
		\begin{cases}
			& (1+|\xi|^2)^\frac{1}{2}|\mathcal{A}_\xi(x,u,\xi)|+|\mathcal{A}(x,u,\xi)|\leq L'(1+|\xi|^2)^\frac{p(x)-1}{2}\\
			
			& \nu' (1+|\xi|^2)^\frac{p(x)-1}{2}|\eta|^2\leq \Big<\mathcal{A}_\xi(x,u,\xi)\eta,\eta\Big>
		\end{cases}
	\end{equation}
	for all $\xi,\eta\in \mathbb{R}^n,u\in\mathbb{R}$ and a.e. $x\in\Omega$ and some positive constants $\nu'\leq L' <+\infty$. $\mathcal{A}_\xi$ denotes differentiation of $\mathcal{A}$ with respect to variable $\xi$. It is clear that structure (\ref{eqh1}) implies (\ref{eqa1}), (\ref{eqa2}) and (\ref{eqa3}).That case $\nu,l,L$ are the function of $n,L',\nu',\gamma_1$ and $\gamma_2$. So it is enough to assume (\ref{eqh1}) instead of (\ref{eqa1},\ref{eqa2}) and (\ref{eqa3}).
	\vskip 1mm In the present paper , it is easy to show that under $p(x)$- growth and coercivity condition ,the natural space for existence of minimizers of the equation is the variable exponent Sobolev space $W^{1,p(.)}(\Omega)$. Such space are relevant in the study of non-Newtonian fluids. For more details about this interesting space,  see \cite{M1,KR2,ER3,ER4,R5,FZ6,RR7,E8,D9,DHHR10}. Apart from these, variable Lebesgue space  has also been used to model  like the thermistor problem \cite{Z13}, quasi-Newtonian fluids \cite{Z65}, magnetostatics \cite{CKMA66} etc. \vskip 1mm
	About last three decades, A good number of Mathematicians show interest on the functionals with $p(x)$-growth for its applicability in the various field like Mathematics, Physics, Technical Studies etc. For example in the modelling of so called electrorheological fluids (see \cite{R11,R12,RR7,MJX21}),modelling of anisotropic materials (see \cite{Z13}), Image Restoring (see \cite{LSC14,CLR15_103}).\vskip 1mm
	In Partial Differential Equations and Calculus of Variations, the topic regualrity is an important and classical segment. The first aim of this paper is to show local Calder\'{o}n-Zygmund estimate of the solution of the given problem (\ref{eq1.1}). This type of estimate had been studied for several years in the case of elliptic and parabolic equations. First this type of work was done in the fundamental paper of T. Iwaniec \cite{I16_62}, where author consider $p$-Laplacean equation whereas in the paper of E.DiBenedetto and J.Manfredi \cite{DM17_55p} had shown Calder\'{o}n-Zygmund estimate for quasilinear elliptic system. Apart from these for constant growth condition see \cite{CP18_60_15_1,DMS19_61_15_2,G20_15_1_24}. In \cite{BDM22_65_15_6} B$\ddot{o}$gerlein, Duzaar, Mingione show optimal Calder\'on-Zygmund theory for $p$-Laplacean type elliptic and parabolic equations under irregular obstacle and this result is extended to the global in \cite{BCW23_66_15_7}. First result of Calder\'on-Zygmund estimate under $p(x)$-growth was proved by Acerbi and Mingione \cite{AM24_63_15_4}. Priori this paper Diening,R$\check{u}$$\check{z}$i$\check{c}$ka \cite{DR25_54p_12} established Calder\'on-Zygmund type estimates for singular integrals in the $L^{p(x)}$ space. later Eleuteri, Habermann \cite{EH26_70_15_11} and Erhardt \cite{E27_71_15_12} proved local Calder\'on-Zygmund estimate for elliptic and parabolic equation respectively. In \cite{BL52}, Byun and Lee proved Calder\'on-Zygmund estimate for nonlinear elliptic double phase problem with $(p(x),q(x))$ - growth conditions.  For global  Calder\'on-Zygmund estimate , Byun , Ok and Ryu \cite{BOR29_64_15_5} consider $p(x)$-Laplacean type whereas in \cite{BCO28_72_15_13} Byun, Cho and Ok proves Global Calder\'on-Zygmund for non-linear elliptic obstacle problems. In \cite{EH26_70_15_11}, authors prove the result with considering "Maximal Function Operator" {\em (Hardy-Littlewood Maximal Function Operator)}. Here we use the technique of \cite{BOR29_64_15_5,BCO28_72_15_13} to prove Calder\'on-Zygmund estimate result. Main novelty of the  paper compare to \cite{BCO28_72_15_13,BOR29_64_15_5} is that we consider obstacle problem depending upon $u$ as well whereas the main difference of this paper with \cite{G20_15_1_24} is that the result dealing with $p(x)$ growth condition where in \cite{G20_15_1_24} author proves result with constant $p$ growth condition.\vskip 1mm
	In recent years Higher Differentiability or Fractional Higher Differentiability results have been attracting a lot of attention starting with some classical papers \cite{N30,N31,N32,KM33,G34,GN35}. Result concerning Higher Differentiability under $p-q$ growth in \cite{CF38_18_1_6,G36_16_1_16,G37_17_1_17}, for subquadratic growth see \cite{G39_21_1_19,N40_26,GG41}. In \cite{GG41} authors prove the result for $p$-growth where in \cite{N40_26} author shows for $p,q$-growth. Apart from these, for standard growth condition see \cite{EN42_13_1_13,EN43_14_1_14,G20_15_1_24}, for Fractional Higher Differentiability result see \cite{CF38_18_1_6,GI44_19_1_25} and for the case of nearly linear growth see \cite{G45_20_1_18}. Most recently , higher order factional differentiability result proved for a system modelling the stationary motion of electro-rheological fluids\cite{MJX21}. In \cite{GI46_41}
	authors prove higher fractional differentiability result for double phase obstracle problem where in \cite{EN43_14_1_14} authors prove Higher Differentiability result in Besov scale of the problem (\ref{eq1.1}) where $f(x,u,\xi)$ is non-differentiable with respect to $x,u$ but only H$\ddot{o}$lder continuous. In \cite{EN43_14_1_14,GI46_41} authors prove the result by introducing "Freezed" functional method.  A.G. Grimaldi \cite{G20_15_1_24} proves higher differentisbility result in Besov Scale for weak solution of (\ref{eq1.1}) under constant $p$-growth condition. As per findings, in the case of variable exponent $p(x)$, the first this type of result was proved for unconstrained case by Giannetti and Napoli \cite{GN47_3}, where author assume suitable Orlicz-Sobolev regularity for $p(.)$ instead of taking the log-H$\ddot{o}$lder continuity. In \cite{FG48_1}, authors prove higher differentisblity result for obstacle problem independent of $u$ variable with $p(x)$ growth. Comparing with the previous results \cite{FG48_1,G20_15_1_24}, this is the first such kind of result upto the knowledge where we prove Higher Differentiablity in Besov scale for $p(x)$-growth of the more general obstacle problem  with additional $u$ dependence. \vskip 1mm
	To prove Calder\'on-Zygmund estiamtes, firstly we use comparison method. In this process , instead of assuming the partial map $x\mapsto F(x,u,\xi)$ is continuous, we assumed comparitively minimul assumption that $x\mapsto\frac{A(x,\xi)}{(1+|\xi|^2)^\frac{p(x)-1}{2}}$  has small BMO semi-norm at (see \ref{eq2.6}) ({\em see for references \cite{BCO28_72_15_13,BOR29_64_15_5}}). Here main difficulty is the presence of variable $p(.)$. To over come this we need to compare the weak solution of equation with variable exponent with the  weak solution of equation with constant exponent (see Lemma \ref{L4.3}) ( {\em for references see \cite{BCO28_72_15_13,BOR29_64_15_5}}).\vskip 1mm
		 
	This Calder\'on-Zygmund result plays a pivotal role during the prove of higher differentiability result. To over come difficulties related to obstacle , we use different quotient method involving both solution and obstacle. To tacle the difficulties related to variable exponent, we uses the boundedness property (see \ref{eq2.1}) of $p(.)$ (see reference \cite{FG48_1}).
	\vskip 1mm
	The assumption on gradient of obstacle to be in {\em Morrey Space} (see subsec. \ref{ss3.5}) is essential because by the result in \cite{EH67}, we have that weak solution of (\ref{eq1.1}) will be locally H$\ddot{o}$lder continuous. This H$\ddot{o}$lder continuity result will  be useful during the prove of the  Calder\'on-Zygmund estimate and higher differentiability results.\vskip 1mm
	
	This paper is organized as follows : In section 2 contians some notations and statements of main results, section 3 describes some spaces including spaces with variable exponent, Besov space etc with some preliminery results, in  section 4 , there are some comparison results with the proof of theorem \ref{T2.1} and section 5 devoted to proof of theorem \ref{T2.2}. 
	\section{Notations and statements of Main results}
	At first we will give the definition of local solution of (\ref{eq1.1}) as follows (see \cite{FG48_1})
	\begin{definition}
		A function $u\in\mathcal{K}_\psi(\Omega)$ is called local solution of obstacle problem (\ref{eq1.1}) in  $\mathcal{K}_\psi(\Omega)$ if $f(x,u,Du)\in L^1_{loc}(\Omega)$ and 
		\begin{equation*}
			\int_{supp(u-\phi)}f(x,u,Du)\leq \int_{supp(u-\phi)}f(x,\phi,D\phi)
		\end{equation*}
		for any $\phi\in \mathcal{K}_\psi(\Omega)$ with $supp(u-\phi)\Subset\Omega$.
	\end{definition}
	\;\;\;\;\;we start with an useful notation. $B_r(x)=B(x;r)=\{z\in\mathbb{R}^n : |z-x|<r\}$ will denote the ball with radius $r$ , centered at $x$. Without any confusion we shall denote $B_r(x)$ as $B_r$. Apart from this we shall denote $f_B$ as the integral average of real valued function $f$ in $B\subset\mathbb{R}^n$ , i.e.
	\begin{equation*}
		(f)_B=f_B=-\!\!\!\!\!\!\int_B f(x)dx=\frac{1}{|B|}\int_B f(x)dx.
	\end{equation*} 
	where $|B|$ is Lebesgue measure of domain $B$.\vskip 1mm
	I assume in this paper that let $p(x):\mathbb{R}^n\rightarrow (1,\infty)$ such that it is a continuous function satisfying 
	\begin{equation}\label{eq2.1}
		1<\gamma_1 :=\inf_\Omega p(x)\leq p(x)\leq \gamma_2 :=\sup_\Omega p(x)<n<+\infty
	\end{equation}
	for some $\gamma_1,\gamma_2$ and for every $x\in \Omega$. Now assume an important condition on $p(x)$ that $p(x)$ is {\em log-H$\ddot{o}$lder} continuous (first introduced by Zhikov\cite{Z13}) in $\Omega$. Then there exist non-decreasing function $\omega: [0,+\infty)\rightarrow[0,+\infty)$, called {\em modulus of continuity} of $p(.)$ such that 
	\begin{equation}\label{eq2.2}
		|p(x)-p(y)|\leq \omega(|x-y|)
	\end{equation}
	for some $x,y\in \mathbb{R}^n$ where $\omega$ satisfying 
	\begin{equation}\label{eq2.3}
		\omega(0)=0 \;\;\;,\;\;\;\lim_{r\rightarrow 0}\omega(r)\log\frac{1}{r}=0
	\end{equation}
	This {\em log-H$\ddot{o}$lder} continuity assumption on $p(.)$ is obvious in the sense that in \cite{Z63} Zhikov shows a counter example where the Lavrentiev's phenomenon may arise which in general can prevent from having regularity of the minimizers.\\[2mm] From now we will assume throughout the paper that $R\leq\frac{1}{2}$. According to \cite[Remark 2]{AM24_63_15_4} we assume that 
	\begin{equation}\label{eq2.4}
		\sup_{0<r\leq R}\omega(r)\log\frac{1}{r}\leq \delta
	\end{equation}
	for some $r<1$ and sufficiently small positive $\delta$ (will be determined later). Now we will take another assumption which is useful to prove of Higher Differentiability result. Assume $\kappa\in L^n\log^nL(\Omega)$ such that 
	\begin{equation}\tag{A5}\label{eqa5}
		|\mathcal{A}(x,u,\xi)-\mathcal{A}(y,v,\xi)|\leq \Bigg(\Big(\kappa(x)+\kappa(y)\Big)|x-y|+\frac{|u-v|^{-\alpha}}{e+|\xi|}\Bigg)(1+|\xi|^2)^\frac{p(x)-1}{2}\log(e+|\xi|^2)
	\end{equation}
\begin{equation}\tag{A6}\label{eqa6}
	l^{-1}(1+|\xi|^2)^\frac{p(x)-1}{2}\leq |\mathcal{A}(x.u,\xi)|
\end{equation}
	for a.e $x,y\in \Omega$ and $u,v\in\mathbb{R},\xi\in\mathbb{R}^n$. $e$ is Euler Constant and $\alpha$ is in Theorem \ref{T3.2}.\vskip 1mm
	Now define 
	\begin{equation}\label{eq2.5}
		A(x,\xi) := \mathcal{A}(x,u(x),\xi)
	\end{equation}
	we assume that $A$ satisfies 
	\begin{equation}\label{eq2.6}
		\sup_{0<r\leq R}\sup_{y\in\mathbb{R}^n}-\!\!\!\!\!\!\int_{B_r(y)}\Theta(A,B_r(y))(x)dx\leq \delta
	\end{equation}
	where 
	\begin{equation}\label{eq2.7}
		\Theta (A,B_r(y))(x)=\sup_{\xi\in\mathbb{R}^n}\Bigg|\frac{A(x,\xi)}{(1+|\xi|^2)^\frac{p(x)-1}{2}}-\Bigg(\frac{A(.,\xi)}{(1+|\xi|^2)^\frac{p(.)-1}{2}}\Bigg)_{B_r(y)}\Bigg|
	\end{equation}
	Now local gradient estimate result as follows,
	\begin{thm}\label{T2.1}
		Let $u\in \mathcal{K}_\psi(\Omega)$ be a weak solution satisfying variational inequalities (\ref{eq1.3}). Let  assumptions $(\ref{eqh1}),(\ref{eqa4}),(\ref{eqa6})$ and $p(.)$ and $\mathcal{A}$ satisfies (\ref{eq2.4}) and (\ref{eq2.6}) respectively for sufficiently small $\delta=\delta(n,\nu',L',\gamma_1,\gamma_2,q)$. Moreover we assume for $q\leq n$, $D\psi\in L^{\gamma_2\tilde{q},\tilde{\lambda}}_{loc}(\Omega)$ for some $\tilde{q}>1$ and $n-\gamma_1<\tilde{\lambda}<n$. Then for $q\in(1,\infty)$, $a^\frac{p(x)}{p(x)-1},|D\psi|^{p(x)}\in L^q_{loc}(\Omega)$, then I have $u\in W^{1,p(x)q}_{loc}(\Omega)$. Indeed the following estimates holds
		\begin{equation}\label{eq2.8}
			\int_{B_R}|Du|^{p(x)q}dx\leq c\int_{B_R}\Big(|Du|^{p(x)}+|u|^{p(x)q}+a^\frac{p(x)q}{p(x)-1}+|D\psi|^{p(x)q}+1\Big)dx
		\end{equation} .
		for $c=c(n,\nu',L',\gamma_1,\gamma_2,q)$
	\end{thm}
	Now the result related to higher differentiability...
	\begin{thm}\label{T2.2}
		Let $\gamma_1>2$,assumptions $(\ref{eqh1}),(\ref{eqa4}),(\ref{eqa5}),(\ref{eqa6})$ holds and $p(.)$ and $\mathcal{A}$ satisfies (\ref{eq2.4}) and (\ref{eq2.6}) respectively for sufficiently small $\delta=\delta(n,\nu',L',\gamma_1,\gamma_2,q)$. For $q>\max\Big\{\frac{3}{2},\frac{\gamma_1}{\gamma_1-2},\frac{n}{n-\frac{\gamma_1}{\gamma_1-2}}\Big\}$, $a^\frac{p(x)}{p(x)-1},|D\psi|^{p(x)}\in L^q_{loc}(\Omega)$. Moreover we assume for $q\leq n$, $D\psi\in L^{\gamma_2\tilde{q},\tilde{\lambda}}_{loc}(\Omega)$ for some $\tilde{q}>1$ and $n-\gamma_1<\tilde{\lambda}<n$, then for weak solution $u\in\mathcal{K}_\psi(\Omega)$ satisfying variational inequalities (\ref{eq1.3}), the following holds
		\begin{equation}\label{eq2.9}
			D\psi\in W^{1,\gamma_2}_{loc}(\Omega)\implies \Big(1+|Du|^2\Big)^\frac{\gamma_1-2}{4}|Du|\in B^\theta_{2,q,loc}(\Omega)
		\end{equation}
		where $\theta<\frac{1}{\gamma_1}$.
	\end{thm}
	\section{Notations, Some Spaces and Preliminery results}
	Let introduce an auxiliary function as following 
	\begin{equation}\label{eq3.1}
		V_p(\eta) :=\Big(1+|\eta|^2\Big)^\frac{p-2}{4}\eta
	\end{equation}
	defined for all $\eta\in\mathbb{R}^n$.\vskip 1mm
	In this paper $c$ denotes a positive coonstant which could be changed in every section , even in every line. Somewhere if needed I use in brackets the dependence of $c$ on relevant parameters.\vskip 1mm
	Throughout this paper for $1\leq p\leq +\infty$, $p'$ will denote the {\em H$\ddot{o}$lder conjugate exponent} of $p$. 
	\vskip 1mm
	\subsection{Space with variable exponent} \label{s3.1}
	\begin{definition}
		Let $\Omega(\subset\mathbb{R}^n)$ be bounded and $p(.):\Omega\rightarrow (1,+\infty)$ be a continuous function. Then generalized Lebesgue space $L^{p(.)}(\Omega)$ or $L^{p(.)}(\Omega,\mathbb{R})$ is defined by 
		\begin{equation*}
			L^{p(.)}(\Omega) := \Bigg\{f:\Omega\rightarrow\mathbb{R} : f\;is\;measurable\;with\;\int_\Omega|f|^{p(x)}dx<+\infty\Bigg\}
		\end{equation*}
		euipped with {\bf Luxemburg Norm} 
		\begin{equation*}
			||f||_{L^{p(.)}(\Omega)} :=\inf \Bigg\{\lambda>o : \int_\Omega\Big|\frac{f(x)}{\lambda}\Big|^{p(x)}dx
			\leq 1 \Bigg\}
		\end{equation*}
		Under this norm $L^{p(.)}(\Omega)$ becomes a separeable Banach space.
	\end{definition}
	If $p(.)$ is unbounded , then variable Lebesgue Space $L^{p(.)}(\Omega)$ is no longer separable space. Moreover class of bounded function of compact support is not dense \cite[p. 4]{CF64}. If no geometric condition on $p(.)$ is imposed than there are only very few properties for the space $L^{p(.)}(\Omega)$. Now we are going to provide an important definition \cite[p. 17]{CF64}
	\begin{definition}
		({\bf The Modular}) Given any domain $\Omega$, $p(.)\in\mathcal{P}(\Omega)$ (where $\mathcal{P}(\Omega)$ is set of Lebesgue measurable functions on $\Omega$) and any Lebesgue measurable function $f$, define the modular functional (or simple the modular) associated with $p(.)$ by
		\begin{equation*}
			\rho(f) = \rho_{p(.),\Omega}(f)=\int_{\Omega\setminus\Omega_\infty}|f(x)|^{p(x)}dx+||f||_{L^\infty(\Omega_\infty)}
		\end{equation*}
	where $\Omega_\infty=\{x\in\Omega : p(x)=\infty\}$
	\end{definition}\vskip 2mm
If $\sup_{\Omega\setminus\Omega_\infty}p(x)<+\infty$, then $L^{p(.)}(\Omega)$ conuncides with the set of functions such that $\rho(f)<\infty$ is finite (see \cite[p. 19]{CF64}). Now let discuss a result (see \cite[prop. 2.41]{CF64})
\begin{prop}\label{p3.1}
	Given $\Omega$ and $p(.)\in\mathcal{P}(\Omega)$. If $f\in L^{p(.)}(\Omega)$, then $f$ is locally integrable.
\end{prop}
	\vskip 1mm
	Now we are going to recall generalized Sobolev spaces $W^{1,p(.)}(\Omega)$ or $W^{1,p(.)}(\Omega,\mathbb{R})$
	\begin{definition}
		The space $W^{1,p(.)}(\Omega)$ is defined as 
		\begin{equation*}
			W^{1,p(.)}(\Omega) := \Bigg\{f\in L^{p(.)}(\Omega) : Df\in L^{p(.)}(\Omega)\Bigg\}
		\end{equation*}
		$Df$ denotes gradient of $f$, equipped with a norm
		\begin{equation*}
			||f||_{W^{1,p(.)}(\Omega)}
			:= ||f||_{L^{p(.)}(\Omega)} + ||Df||_{L^{p(.)}(\Omega)}
		\end{equation*}
		space $W^{1,p(.)}(\Omega)$ with this norm is Banach space.
	\end{definition}
	\vskip 1mm
	For more details about these spaces see \cite{DHHR10, CF64}.
	\vskip 5mm
	For the Auxiliary function $V_p$, let recall the following result (see proof in \cite[Lemma 8.3]{G51}) 
	\begin{lem}\label{L3.1}
		Let $1<p<+\infty$. There exists a constant $C=C(n,p)>0$ such that 
		\begin{equation*}
			C^{-1}\Big(1+|\xi|^2+|\eta|^2\Big)^\frac{p-2}{2}\leq \frac{|V_p(\xi)-V_p(\eta)|^2}{|\xi-\eta|^2}\leq C\Big(1+|\xi|^2+|\eta|^2\Big)^\frac{p-2}{2}
		\end{equation*}
		for any $\xi,\eta\in\mathbb{R}^n$ , $\xi\neq\eta$.
	\end{lem}
	Additionally I will often use the inequality during the proof of 
	\begin{lem}\label{L3.2} (see \cite{GN47_3,FG48_1})
		For every $s,t>0$ and for every $\alpha,\epsilon,\gamma>0$, I have 
		\begin{equation*}
			st\leq \epsilon s\log^\alpha(e+\gamma s)+\frac{t}{\gamma}\Bigg[\exp\Big(\frac{t}{\epsilon}\Big)^\frac{1}{\alpha}-1\Bigg]
		\end{equation*}
	\end{lem}
	\subsection{Difference quotient}
	for every function $f:\mathbb{R}^n\rightarrow\mathbb{R}$, the finite difference operator defined with respect to $x_s$ by 
	\begin{equation*}
		\tau_{s,h} = f(x+he_s)-f(x)
	\end{equation*}
	where $h\in \mathbb{R}^n$ and $e_s$ is the unit vector denoting direction of $x_s$ axis and $s\in\{1,2,3,...,n\}$.\vskip 1mm
	Some properties of the finite difference quotient will be recalled for future purpose. Let start with some elementary properties which can be found , see for example in \cite{G51,G20_15_1_24}.
	\begin{prop}
		Let two functions $f$ and $g$ be such that $f,g\in W^{1,p}(\Omega)$, with $p\geq 1$, and let consider the set 
		\begin{equation*}
			\Omega_{|h|}=\{x\in\Omega : |h|<dist(x,\partial\Omega)\}.
		\end{equation*}
		Then,\\[0.8mm]
		(i) $\tau_hf\in W^{1,p}(\Omega_{|h|})$ and 
		\begin{equation*}
			D_i(\tau_hf)=\tau_h(D_if)
		\end{equation*}
		(ii) If at least one of the functions $f$ and $g$ has support contained in $\Omega_{|h|}$, then 
		\begin{equation*}
			\int_\Omega f\tau_h g dx=\int_\Omega g\tau_h f dx
		\end{equation*}
		(iii) I have 
		\begin{equation*}
			\tau_h(fg)(x)=f(x+h)\tau_h g(x)+g(x)\tau_h f(x)
		\end{equation*}
	\end{prop}
	The the following result regarding finite difference operator is a kind of integral version of Lagrange Theorem
	\begin{lem}\label{L3.3}
		If $0<\rho<R,|h|<\frac{R-\rho}{2},1<p<+\infty$ and $f\in W^{1,p}(B_R)$, then
		\begin{equation*}
			\int_{B_\rho}|\tau_hf(x)|^pdx\leq c(n,p)|h|^p\int_{B_R}|Df(x)|^pdx
		\end{equation*}
		Moreover,
		\begin{equation*}
			\int_{B_\rho}|f(x+h)|^pdx\leq\int_{B_R}|f(x)|^pdx
		\end{equation*}
	\end{lem}
	For the proof of first part see for example in \cite[Section 5.8.2 , Theorem 3(i)]{E53_127}.\vskip 1mm
	Now, recall the fundamental Sobolev embedding property.
	\begin{lem}
		Let $f:\mathbb{R}^n\rightarrow\mathbb{R}^N, f\in L^p(B_R)$ with $1<p<n$. Suppose that there exist $\rho\in (0,R)$ and $M>0$ such that 
		\begin{equation*}
			\sum_{s=1}^{n}\int_{B_\rho}|\tau_{s,h}f(x)|^pdx\leq M^p|h|^p
		\end{equation*}
		for every $h$ with $|h|<\frac{R-\rho}{2}$. Then $f\in W^{1,p}(B_\rho)\cap L^\frac{np}{n-p}(B_\rho)$. Moreover,
		\begin{equation*}
			||Df||_{L^p(B_\rho)}\leq M
		\end{equation*}
		and 
		\begin{equation*}
			||f||_{L^\frac{np}{n-p}(B_\rho)}\leq c(M+||Df||_{L^p(B_R)})
		\end{equation*}
		with $c=c(n,N,p,\rho,R)$
	\end{lem}
	For the proof see for example \cite[Lemma 8.2]{G51}.\vskip 1mm
	Now conclude this subsection by recalling a fractional version of previous lemma (see \cite{KM33})
	\begin{lem}
		Let $f\in L^2(B_R)$. Suppose that there exist $\rho\in (0,R),0<\alpha<1$ and $M>0$ such that
		\begin{equation*}
			\sum_{s=1}^{n}\int_{B_\rho}|\tau_{s,h}f(x)|^pdx\leq M^2|h|^{2\alpha}
		\end{equation*} 
		or every $h$ with $|h|<\frac{R-\rho}{2}$. Then $f\in L^\frac{2n}{n-2\beta}(B_\rho)$ for every $\beta\in (0,\alpha)$ and 
		\begin{equation*}
			||f||_{L^\frac{2n}{n-2\beta}(B_\rho)}\leq c(M+||f||_{L^2(B_R)})
		\end{equation*}
		with $c=c(n,N,R,\rho,\alpha,\beta)$.
	\end{lem}
	\subsection{Spaces $L^p\log^\sigma L$}
	\begin{definition}
		For all $p\geq 1$ and $\sigma\in\mathbb{R}$ , the space $L^p\log^\sigma L(\Omega)$ is defined by 
		\begin{equation*}
			L^p\log^\sigma L := \Big\{f:\Omega\rightarrow\mathbb{R}: f\;is\;measurable\;and\;\int_\Omega|f|^p\log^\sigma(e+|f|)dx<+\infty\Big\}
		\end{equation*}
		and it becomes Banach space under Luxemburg norm
		\begin{equation}
			||f||_{L^p\log^\sigma L} :=\inf \Bigg\{\lambda>o : \int_\Omega\Big|\frac{f(x)}{\lambda}\Big|^p\log^\sigma \Big(e+\frac{|f|}{\lambda}\Big)dx
			\leq 1 \Bigg\}
		\end{equation}
	\end{definition}\vskip 2mm 
	Now we want to recall some properties and lemmas regarding space $L\log^\sigma L$ (see and corresponding references in  \cite{AM24_63_15_4}).
	\begin{lem}
		For $p>1$ and $\sigma>1$, then if $f\in L^p(\Omega)$ then $f\in L\log^\sigma L$, that is there exist a constant $c=c(p)>0$ such that 
		\begin{equation}
			-\!\!\!\!\!\!\int_\Omega |f|\log^\sigma \Bigg(e+\frac{|f|}{-\!\!\!\!\!\int_\Omega|f|dx}\Bigg)dx\leq c\Bigg(-\!\!\!\!\!\!\int_\Omega|f|^pdx\Bigg)^\frac{1}{p}
		\end{equation} 
	\end{lem}\vskip 3mm
	Before proceeding to lemmas let recall , due to  \cite{I57,IV58} 
	\begin{equation}
		[f]_{L\log^\sigma L} := -\!\!\!\!\!\!\int_\Omega |f|\log^\sigma \Bigg(e+\frac{|f|}{-\!\!\!\!\!\int_\Omega|f|dx}\Bigg)dx
	\end{equation}
	\begin{lem}
		For every $f\in L\log^\sigma L$, then  both norm $[f]_{L\log^\sigma L}$ and $||f||_{L\log^\sigma L}$ are comparable. That is there exists $c=c(\beta)\geq 1$ , independent of $\Omega$ and $f$ such that 
		\begin{equation}
			\frac{1}{c}||f||_{L\log^\sigma L}\leq [f]_{L\log^\sigma L}\leq c ||f||_{L\log^\sigma L}
		\end{equation} 
	\end{lem}
	
	\subsection{Besov-Lipscchitz spaces}
	Let $v:\mathbb{R}^n\rightarrow\mathbb{R}$ be a function and $h\in\mathbb{R}^n$. As in \cite{H60}, given $0<\alpha<1$ and $1\leq p,q<\infty$, we say that $v$ belogs to the Besov space $B^\alpha_{p,q}(\mathbb{R}^n)$ if $v\in L^p(\mathbb{R}^n)$ and 
	\begin{equation}
		||v||_{B^\alpha_{p,q}(\mathbb{R}^n)} = ||v||_{L^p(\mathbb{R}^n)}+[v]_{B^\alpha_{p,q}(\mathbb{R}^n)}<\infty
	\end{equation}
	where 
	\begin{equation}
		[v]_{B^\alpha_{p,q}(\mathbb{R}^n)} = \Bigg(\int_{\mathbb{R}^n}\Bigg(\int_{\mathbb{R}^n}\frac{|v(x+h)-v(x)|^p}{|h|^{\alpha p}}dx\Bigg)^\frac{q}{p}\frac{dh}{|h|^n}\Bigg)^\frac{1}{q}<\infty
	\end{equation}
	Equivalently ,simply say that $v\in L^p(\mathbb{R}^n)$ and $\frac{\tau_hv}{|h|^\alpha}\in L^q\Big(\frac{dh}{h^n};L^p(\mathbb{R}^n)\Big)$. As ususal, if one simply integrates for $h\in B(0,\delta)\;,\;\delta>0$, $\delta$ is fixed then an equivalent norm is obtained , because
	\begin{equation}
		\Bigg(\int_{|h|\geq \delta}\Bigg(\int_{\mathbb{R}^n}\frac{|v(x+h)-v(x)|^p}{|h|^{\alpha p}}dx\Bigg)^\frac{q}{p}\frac{dh}{|h|^n}\Bigg)^\frac{1}{q}\leq c(n,\alpha,p,q,\delta)||v||_{L^p(\mathbb{R}^n)}
	\end{equation}
	Similarly , we can say that $v\in{B^\alpha_{p,q}(\mathbb{R}^n)}$ if $v\in L^p(\mathbb{R}^n)$ and 
	\begin{equation}
		[v]_{B^\alpha_{p,\infty}(\mathbb{R}^n)}=\sup_{h\in\mathbb{R}^n}\Bigg(\int_{\mathbb{R}^n}\frac{|v(x+h)-v(x)|^p}{|h|^{\alpha p}}dx\Bigg)^\frac{1}{p}<\infty
	\end{equation}
	Again, simply one can take supremum over $|h|\leq \delta$ and obtain an equivalent norm. By construction, $B^\alpha_{p,q}(\mathbb{R}^n)\subset L^p(\mathbb{R}^n)$. One also has the following version of Sobolev embeddings (see proof in \cite{H60})
	\begin{lem}
		Suppose that $0<\alpha<1$, then if $1<p<\frac{n}{\alpha}$ and $1\leq q \leq p^*_\alpha=\frac{np}{n-\alpha p}$, then there is a continuous embedding $B^\alpha_{p,q}(\mathbb{R}^n)\subset L^{p^*_\alpha}(\mathbb{R}^n)$.
	\end{lem}\vskip 1mm
	The proof of following lemma can be found in \cite{BCGON61}.
	\begin{lem}\label{L3.9}
		A function $v\in L^p_{loc}(\Omega)$ belongs to the local Besov space $B^\alpha_{p,q,loc}$ iff
		\begin{equation}
			\Bigg|\Bigg|\frac{\tau_hv}{|h|^\alpha}\Bigg|\Bigg|_{L^q\big(\frac{dh}{|h|^n};L^p(B)\big)}<\infty
		\end{equation}
		for any ball $ B\subset 2B\Subset \Omega$ with radius $r_B$. Here the measure $\frac{dh}{|h|^n}$ is restricted to the ball $B(0,r_B)$ on the $h$-space.
	\end{lem}
	Besov-Lipschitz regualrity of a function can be characterized in pointwise terms. Given measurable function $v:\in\mathbb{R}^n\rightarrow\mathbb{R}$, a {\em fractional $\alpha$-Hajlasz gradient for $v$} is a sequence $\{g_k\}_k$ of measurable , non-negative functions $g_k:\mathbb{R}^n\rightarrow\mathbb{R}$, together with a null set $N\subset\mathbb{R}^n$, such that the inequality 
	\begin{equation}
		|v(x)-v(y)|\leq (g_k(x)+g_k(y))|x-y|^\alpha
	\end{equation}
	holds whenever $k\in\mathbb{Z}$ and $x,y\in\mathbb{R}^n\setminus N$ are such that $2^{-k}\leq |x-y|\leq 2^{-k+1}$. This $\{g_k\}_k\in l^q(\mathbb{Z};L^p(\mathbb{R}^n))$ if 
	\begin{equation}
		||\{g_k\}_k||_{l^q(L^p)}=\Bigg(\sum_{k\in\mathbb{Z}}||g_k||^q_{L^p(\mathbb{R}^n)}\Bigg)<\infty
	\end{equation}
	Now a theorem proved in \cite{KYZ61} as follows
	\begin{thm}
		Let $0<\alpha<1,1\leq p<\infty$ and $1\leq q\leq \infty$. Let $v\in L^p(\mathbb{R}^n)$. One has $v\in B^\alpha_{p,q}(\mathbb{R}^n)$ iff there exists  a fractional $\alpha$-Hajlasz gradient $\{g_k\}_k\in l^q(\mathbb{Z} ; L^p(\mathbb{R}^n))$ for $v$. Moreover 
		\begin{equation}
			||v||_{B^\alpha_{p,q}(\mathbb{R}^n)} \simeq \inf ||\{g_k\}||_{l^q(L^p)}
		\end{equation}
		where the infimum runs over all possible fractional $\alpha$-Hajlasz gradient for $v$.
	\end{thm}
\subsection{Morrey Space}\label{ss3.5}
Let $1\leq p<+\infty$ and $\lambda\geq 0$. Let $\Omega$ is bounded open subset of $\mathbb{R}^n$. By $L^{p,\lambda}(\Omega)$ we denote a linear space of functions $u\in L^p(\Omega)$ such that 
\begin{equation}
	||u||_{L^{p,\lambda}(\Omega)} := \Bigg\{\sup_{x_0\in\Omega, 0<\rho<diam(\Omega)} \rho^{-\lambda}\int_{\Omega(x_0,\rho)}|u(x)|^pdx\Bigg\}^{1/p}<+\infty
\end{equation}
where we set $\Omega(x_0,\rho) := \Omega\cap B(x_0,\rho)$. Finally we can have from \cite[subsec. 2.1]{G20_15_1_24}, $||u||_{L^{p,0}(\Omega)}=||u||_{L^p(\Omega)}$, so that $ L^{p,0}(\Omega)=L^p(\Omega)$. Using H$\ddot{o}$lder inequality, one have that if $s\geq p$ and $\frac{n-\lambda}{p}\geq \frac{n-\kappa}{s}$ the following holds: 
\begin{equation}
	||u||_{L^{p,\lambda}(\Omega)}\leq\; diam(\Omega)^{\frac{n-\lambda}{p}-\frac{n-\kappa}{s}}||u||_{L^{s,\kappa}(\Omega)}
\end{equation}
and therefore the inclusion
\begin{equation}
	L^{s,\kappa}(\Omega)\subset L^{p,\lambda}(\Omega)
\end{equation}
is continuous.\\[2mm]
Now we want to provide a result which we will use intensively in the sequel. The result can be found in \cite[Th. 2.7]{EH67}
\begin{thm}\label{T3.2}
	Let $u\in W^{1,p(.)}(\Omega)$ be alocal minimizer of the functional (\ref{eq1.1}) in the class of (\ref{eq1.2}), where $\psi\in W^{1,1}_{loc}(\Omega)$ is a given obstacle function fulfilling 
	\begin{equation}
		D\psi\in L^{\gamma_2\tilde{q},\tilde{\lambda}}_{loc}(\Omega)
	\end{equation}
with some $\tilde{q}>1$ and $n-\gamma_1<\tilde{\lambda}<n$, where $\gamma_1$ and $\gamma_2$ have been introduced in (\ref{eq2.1}). Suppose moreover that the $\mathcal{A}$ satisfies (\ref{eqa3}) and (\ref{eqa6}) and also the function $p(.)$ fulfills assumptions (\ref{eq2.2}) and (\ref{eq2.4}). Then $u\in \mathcal{C}^{0,\alpha}_{loc}(\Omega)$ for some $\alpha\in (0,1)$.
\end{thm}
Then we have for every $\Omega'\Subset \Omega$ , there exist some constant $D$ such that 
\begin{equation}\label{eq3.18}
	|u(x)-u(y)|\leq D|x-y|^\alpha
\end{equation}
for some $x,y\in\Omega'$.
	\section{Local Gradient Estimates}
	In this section we will prove local Calder\'on-Zygmund estimates for solution of variational inequality (\ref{eq1.3}). Priori this we will discuss some preliminery results (see for references \cite{BCO28_72_15_13,BOR29_64_15_5,G20_15_1_24}).
	\subsection{Comparison estimates by approximation}
	For sake convenience, we define $\Psi = a(x)^\frac{1}{p(x)-1}+|D\psi(x)|+1$. Fix $R\in (0,1)$ in such a way that $B_R\Subset\Omega$. $A$ is defined as in (\ref{eq2.5}).\vskip 1mm
	For the weak solution $u\in W^{1,p(.)}(\Omega)$ to the variational inequality (\ref{eq1.3}) , let $k\in W^{1,p(x)}(B_{4r})$ be the weak solution of 
	\begin{equation}\label{eq4.1}
		\begin{cases}
			div\;A(x,Dk) = div\;A(x,D\psi) \;\;\;\; in\;B_{4r}\\
			k=u \hspace{40mm} on\;\partial B_{4r}
		\end{cases}
	\end{equation}
	and $w\in W^{1,p(.)(\Omega)}$ be the weak solution of 
	\begin{equation}\label{eq4.2}
		\begin{cases}
			div\;A(x,Dw)=0\;\;\;in\;B_{4r}\\
			w=k\hspace{20mm} on\; \partial B_{4r}
		\end{cases}
	\end{equation}
	Now the prove the following (modifying Lemma 3.2 in \cite{G20_15_1_24} based on \cite{BCO28_72_15_13})
	\begin{lem}\label{L4.1}
		Let $\epsilon\in (0,1)$ and $\lambda\geq 1$. There exists $\delta=\delta(n,\nu,L,\gamma_1,\gamma_2,\epsilon)>0$ such that if  such that if
		\begin{equation}\label{eq4.3}
			-\!\!\!\!\!\!\int_{B_{4r}}(|Du|^{p(x)}+|u|^{p(x)})dx \leq \lambda\;\;\;and\;\;\;-\!\!\!\!\!\!\!\int_{B_{4r}} \Psi^{p(x)} dx
			\leq \delta\lambda
		\end{equation}
		{\em then }
		\begin{equation}\label{eq4.4}
			-\!\!\!\!\!\!\int_{B_{4r}}|Dw|^{p(x)}dx\leq c\lambda
		\end{equation} 
		{\em for some constant $c=c(n,\nu,L,\gamma_1,\gamma_2,\epsilon)$ and}
		\begin{equation}\label{eq4.5}
			-\!\!\!\!\!\!\int_{B_{4r}}|Du-Dw|^{p(x)}dx\leq \epsilon\lambda
		\end{equation}
	\end{lem}
	{\bf Proof.} By standard energy estimates (as local solutions continuously depend on data(s)) for eqations (\ref{eq4.1}) and (\ref{eq4.2}) , we have from (\ref{eq4.3}) that 
	\begin{equation}\label{eq4.6}
		-\!\!\!\!\!\!\int_{B_{4r}}|Dk|^{p(x)}dx\leq -\!\!\!\!\!\!\int_{B_{4r}} c\Big(|D\psi|^{p(x)}+|Du|^{p(x)}+1\Big)\leq c\lambda
	\end{equation}
	and so 
	\begin{equation}\label{eq4.7}
		-\!\!\!\!\!\!\int_{B_{4r}}|Dw|^{p(x)} dx\leq c \Big(-\!\!\!\!\!\!\!\int_{B_{4r}}|Dk|^{p(x)}dx+1\Big)\leq c\lambda
	\end{equation}
	In view of \cite[Lemma 4.1]{BCO28_72_15_13} and equations (\ref{eq4.1}) , we will have that $k\geq\psi$ a.e. in $B_{4r}$. we extend $k$ by $u$ in $\Omega-B_{4r}$, hence , we have $k\in \kappa_\psi(\Omega)$. Therefore we  can take $\psi=k$ in (\ref{eq1.3}), thus obtaining
	\begin{equation}\label{eq4.8}
		\int_{B_{4r}}\big<A(x,Du),D(k-u)\big> dx\geq \int_{B_{4r}}\mathcal{B}(x,u,Du)(k-u)dx
	\end{equation}
	Choosing $k-u\in W_0^{1,p(.)}(B_{4r})$ as a test function for (\ref{eq4.1}), we also have 
	\begin{equation}\label{eq4.9}
		\int_{B_{4r}}\big<A(x,Dk),D(k-u)\big>=\int_{B_{4r}}\big<A(x,D\psi),D(k-u)\big>dx
	\end{equation} 
	now subtracting (\ref{eq4.8}) form (\ref{eq4.9}) to find 
	\begin{equation}\label{eq4.10}
		-\!\!\!\!\!\!\int_{B_{4r}}\big<A(x,Dk)-A(x,Du),D(k-u)\big>dx \leq -\!\!\!\!\!\!\int_{B_{4r}} \big<A(x,d\psi),D(k-u)\big> dx - -\!\!\!\!\!\!\int_{B_{4r}}\mathcal{B}(x,u,Du)(k-u)dx
	\end{equation}
	Now estimating  left hand side of (\ref{eq4.10}). 
	Since we already assume $p(x)\geq 2$ then by virtue of (\ref{eqa1}) we shall have that,
	\begin{equation}\label{eq4.11}
		-\!\!\!\!\!\!\int_{B_{4r}}|Dk-Du|^{p(x)}\leq \frac{1}{\tilde{\nu}}-\!\!\!\!\!\!\int_{B_{4r}}\big<A(x,Dk)-A(x,Du),D(k-u)\big>dx
	\end{equation}
	Now next estiamte the right hand side in (\ref{eq4.10}) . By assumption (\ref{eqa3}) and (\ref{eqa4}), from Young's inequality, Sobolev embedding theorem  variable exponent (see \cite{DHHR10}), denoting by $c_0$ the Poincar\'e Constant , we have that for any positive numbers $\eta$ and $\theta$,
	\begin{equation*}
		-\!\!\!\!\!\!\int_{B_{4r}} \big<A(x,D\psi),D(k-u)\big>dx- 	-\!\!\!\!\!\!\int_{B_{4r}} \mathcal{B}(x,u,Du)(k-u)dx\leq c	-\!\!\!\!\!\!\int_{B_{4r}}(1+|D\psi|^{p(x)-1})|Dk-Du|dx
		+ 	-\!\!\!\!\!\!\int_{B_{4r}}(|Du|^{r(x)}+|u|^{r(x)}+a)|k-u|dx
	\end{equation*}
	\begin{equation*}\footnote{use Young's Inequality by assuming $p(x)=p(x)$, $p'(x)=\frac{p(x)}{p(x)-1}$}
		\leq\eta -\!\!\!\!\!\!\int_{B_{4r}}|Dk-Du|^{p(x)}dx+c(\eta)-\!\!\!\!\!\!\int_{B_{4r}}(1+|D\psi|^{p(x)})dx+\frac{\eta}{c_0} -\!\!\!\!\!\!\int_{B_{4r}}|k-u|^{p(x)}dx+c(\eta)-\!\!\!\!\!\int_{B_{4r}}(|Du|^{r(x)p'(x)}+|u|^{r(x)p'(x)}+a^{p'(x)})dx
	\end{equation*}
	\begin{equation*}\footnote{Using Poincar\'e Inequality}
		\leq \eta -\!\!\!\!\!\!\int_{B_{4r}}|Dk-Du|^{p(x)}dx+ c(\eta) -\!\!\!\!\!\!\int_{B_{4r}}(1+|D\psi|^p(x)+a^{p'(x)})dx+\eta -\!\!\!\!\!\!\int_{B_{4r}}|Dk-Du|^{p(x)}dx+c(\eta)-\!\!\!\!\!\int_{B_{4r}}(|Du|^{r(x)p'(x)}+|u|^{r(x)p'(x)})dx
	\end{equation*}
	\begin{equation*}
		\leq 2\eta -\!\!\!\!\!\!\int_{B_{4r}}|Dk-Du|^{p(x)}dx+ c(\eta) -\!\!\!\!\!\!\int_{B_{4r}}(1+|D\psi|^{p(x)}+a^{p'(x)})dx+\theta c(\eta)-\!\!\!\!\!\!\int_{B_{4r}}(|Du|^{p(x)}+|u|^{p(x)})dx + c(\eta)c(\theta)
	\end{equation*}
	\begin{equation}\label{eq4.12}
		\leq2\eta-\!\!\!\!\!\!\int_{B_{4r}}|Dk-Du|^{p(x)}dx+c(\eta,\theta)\delta\lambda+c(\eta)\theta\lambda
	\end{equation}
	Coombining (\ref{eq4.10},\ref{eq4.11},\ref{eq4.12}), we derive that
	\begin{equation}\label{eq4.13}
		-\!\!\!\!\!\!\int_{B_{4r}} |Dk-Du|^{p(x)}d\leq \frac{1}{\tilde{\nu}}\Big(2\eta -\!\!\!\!\!\!\int_{B_{4r}} |Dk-Du|^{p(x)}dx + C(\eta,\theta)\delta\lambda+C(\eta)\theta\lambda\Big)
	\end{equation} 
	Choosing 
	\begin{equation*}
		\eta=\frac{\tilde{\nu}}{4}\;,\;\theta=\frac{\tilde{\nu}\epsilon}{4c(\eta)}\;,\;\delta<\frac{\tilde{\nu}\epsilon}{c(\eta,\theta)}
	\end{equation*}
	we will have 
	\begin{equation}\label{eq4.14}
		-\!\!\!\!\!\!\int_{B_{4r}} |Dk-Du|^{p(x)}dx\leq \epsilon\lambda
	\end{equation}
	on other hand , by taking $w-k\in W^{1,p(.)}_0(B_{4r})$ as a test function for (\ref{eq4.1}) and (\ref{eq4.2}) and then substracting, we have 
	\begin{equation*}
		-\!\!\!\!\!\!\int_{B_{4r}}\big<A(x,Dw)-A(x,Dk),D(w-k)\big>dx=-\!\!\!\!\!\!\int_{B_{4r}}\big<A(x,D\psi),D(k-u)\big>dx
	\end{equation*}
	In a similar way of estimating (\ref{eq4.14}), we can find $\delta=\delta(n,\nu,L,\gamma_1,\gamma_2,\epsilon)>0$ such that 
	\begin{equation*}
		-\!\!\!\!\!\!\int_{B_{4r}} |Dw-Dk|^{p(x)}\leq \epsilon\lambda \hspace{10mm} \Box
	\end{equation*}
	Now taking $k-u\in W^{1,p(x)}_0(B_{4r})$ as a test function for (\ref{eq4.1}) , we have 
	\begin{equation*}
		-\!\!\!\!\!\!\int_{B_{4r}}\big<A(x,Dk),Dk-Du\big>=-\!\!\!\!\!\!\int_{B_{4r}}\big<A(x,D\psi),Dk-Du\big>dx
	\end{equation*}
	and so 
	\begin{equation}\label{eq4.15}
		-\!\!\!\!\!\!\int_{B_{4r}}\big<A(x,Dk)-A(x,du),Dk-Du\big>=-\!\!\!\!\!\!\int_{B_{4r}}\big<A(x,D\psi)-A(x,Du),Dk-Du\big>dx
	\end{equation}
	By the ellipticity assumptions (\ref{eqa1}) and (\ref{eqa2}) ,we get
	\begin{equation*}
		{\tilde{\nu}}-\!\!\!\!\!\!\int_{B_{4r}}|Dk-Du|^{p(x)}dx\leq -\!\!\!\!\!\!\int_{B_{4r}} \big<A(x,D\psi-A(x,Du),Dk-Du)dx
	\end{equation*}
	\begin{equation*}
		\leq -\!\!\!\!\!\!\int_{B_{4r}} |A(x,D\psi)-A(x,Du)||Dk-Du|
	\end{equation*}
	\begin{equation*}
		\leq -\!\!\!\!\!\!\int_{B_{4r}}(1+|D\psi|^{p(x)-2}+|Du|^{p(x)-2})|Dk-Du|^2dx
	\end{equation*}
	\begin{equation*}\footnote{Use Young's Inequality with $p=p(x)/2$ and $p'=\frac{p(x)}{p(x)-2}$}
		\leq\theta  -\!\!\!\!\!\!\int_{B_{4r}}|Dk-Du|^{p(x)}+C(\theta)  -\!\!\!\!\!\!\int_{B_{4r}} (1+|D\psi|^{p(x)}+|Du|^{p(x)})dx
	\end{equation*}
	So choosing $\theta=\frac{{\tilde{\nu}}}{2}$ ,we can reabsorb the first integral in the right hand side by left hand side and we obtain
	\begin{equation*}
		-\!\!\!\!\!\!\int_{B_{4r}} |Dk-Du|^{p(x)}\leq c-\!\!\!\!\!\!\int_{B_{4r}}(1+|D\psi|^{p(x)}+|Du|^{p(x)})dx
	\end{equation*}
	that by virtue of assumption (\ref{eq4.3}) yields
	\begin{equation}\label{eq4.16}
		-\!\!\!\!\!\!\int_{B_{4r}}|Dk|^{p(x)}\leq c \big(-\!\!\!\!\!\!\int_{B_{4r}}|D\psi|^{p(x)}+-\!\!\!\!\!\!\int_{B_{4r}}|Du|^{p(x)}+1\big)\leq c\lambda
	\end{equation}
	On the other hand , taking $w-k\in W^{1,p(x)}_0(B_{4r})$ as test function in (\ref{eq4.2}), I have 
	\begin{equation*}
		-\!\!\!\!\!\!\int_{B_{4r}} \big<A(x,Dw),Dw-Dk\big>dx=0
	\end{equation*}
	and so 
	\begin{equation*}
		-\!\!\!\!\!\!\int_{B_{4r}} \big<A(x,Dw)-A(x,Dk),Dw-Dk\big>dx=--\!\!\!\!\!\!\int_{B_{4r}}\big<A(x,Dk),Dw-Dk\big>dx
	\end{equation*}
	By the use of ellipticity assumption (\ref{eqa1}) in the left hand side and assumption (\ref{eqa3}) in the right hand side, I will have 
	\begin{equation*}
		\tilde{\nu} -\!\!\!\!\!\!\int_{B_{4r}} |Dw-Dk|^{p(x)}dx\leq -\!\!\!\!\!\!\int_{B_{4r}} \big<A(x,Dw)-A(x,Dk),Dw-Dk\big>dx
	\end{equation*} 
	\begin{equation*}
		\leq -\!\!\!\!\!\!\int_{B_{4r}} (1+|Dk|^{p(x)-1})|Dw-Dk|dx
	\end{equation*}
	\begin{equation*}\footnote{Use Young's Inequality with $p=p(x)$ and $p'=\frac{p(x)}{p(x)-1}$}
		\leq \theta -\!\!\!\!\!\!\int_{B_{4r}}|Dw-Dk|^{p(x)}dx+C(\theta) -\!\!\!\!\!\!\int_{B_{4r}}(1+|Dk|^{p(x)})dx
	\end{equation*}
	So, choosing $\theta=\frac{\tilde{\nu}}{2}$, we obtain
	\begin{equation*}
		-\!\!\!\!\!\!\int_{B_{4r}}|Dw-Dk|^{p(x)}dx\leq c -\!\!\!\!\!\!\int_{B_{4r}}(1+|Dk|^{p(x)})dx
	\end{equation*}
	that by (\ref{eq4.16}) yields
	\begin{equation*}
		-\!\!\!\!\!\!\int_{B_{4r}}|Dw|^{p(x)}dx\leq c\big(-\!\!\!\!\!\!\int_{B_{4r}}|Dk|^{p(x)}dx+1\big)dx\leq c\lambda\hspace{10mm}\Box
	\end{equation*}\\[3mm]
	Let discuss a Lemma as follows...
	\begin{lem}\label{L4.2}
		Let assume that (\ref{eqa1},\ref{eqa3},\ref{eq2.4},\ref{eq2.1}) For any fixed $r>0$ satisfying 
		\begin{equation*}
			r\leq \frac{R}{4}\;\;\; and \;\;\; \omega(4r)\leq \sqrt{\frac{n+1}{n}}-1
		\end{equation*}
		Then there exists $c\equiv c(n,\nu',L',\gamma_1,\gamma_2)$ and $c_g\equiv c_g(n,\nu',L',\gamma_1,\gamma_2)$ such that for any 
		\begin{equation*}
			0<\sigma\leq \min\Bigg\{\frac{c_g}{M_1^{4\omega(4r)/{\gamma_1}}},1\Bigg\}=:\sigma_0\;\;\;\;with\;\;M_1 := \int_{B_{4r}}|Dw|^{p(x)}+1
		\end{equation*}
		Then for every $B_\rho\subseteq B_{4r}\subset\subset \Omega$, it holds 
		\begin{equation}\label{eq4.17}
			-\!\!\!\!\!\!\int_{B_\rho}|Dw|^{p(x)(1+\sigma)}dx\leq c \Biggl\{\Bigg(-\!\!\!\!\!\!\!\int_{B_2\rho}|Dw|^{p(x)}dx\Bigg)^{(1+\sigma)}+1\Biggr\}
		\end{equation}
	\end{lem}
	{\bf Proof.} see proof in \cite[Theorem 5]{AM24_63_15_4}.\vskip 1mm
	If $u\in W^{1,p(.)}(\Omega)$ is a weak solution of the variational problem (\ref{eq1.1}) under assumed conditins on $\psi,\mathcal{A},\mathcal{K}_\psi$, then one can have the following estimate ({\em solution continuously depends on data }) 
	\begin{equation}\label{eq4.18}
		\int_\Omega|Du|^{p(x)}\leq c\int_\Omega\big(|D\psi|^{p(x)}+1\big)dx
	\end{equation}
	Now we observe that using (\ref{eq4.18}) as following,
	\begin{equation*}
		\int_{B_{4r}}|Dw|^{p(x)}dx+1\leq c\int_{B_{4r}}(|Dk|^{p(x)}+1)dx+1
	\end{equation*}
	\begin{equation*}
		\hspace{49mm}\leq c\int_{B_{4r}}(|Du|^{p(x)}+|D\psi|^{p(x)}+1)dx +1
	\end{equation*}
	\begin{equation*}
		\hspace{39mm}\leq c_*\int_{B_{4r}}(|D\psi|^{p(x)}+1)dx+1
	\end{equation*}
	\begin{equation*}
		\hspace{48mm}\leq c_* \int_{B_{4r}}(|a|^{\frac{1}{p(x)-1}}+|D\psi|+1)^{p(x)}dx+1
	\end{equation*}
	\begin{equation}\label{eq4.19}
		\hspace{48mm}\leq c_*M\;\;\;\;\;where \;\;M=\int_{B_{4r}}|\Psi|^{p(x)}dx+1
	\end{equation}
	Let write 
	\begin{equation*}
		p^-=\inf_{x\in B_{4r}} p(x)\;\;and\;\;p^+=\sup_{x\in B_{4r}} p(x)
	\end{equation*}
	Now let us provide an estimate using the Lemma \ref{L4.2} and it will provide the relation between variable exponent $p(.)$ and constant exponent $p^+$ as follows\vskip 0.8mm
	{\em For $w\in W^{1,p^+}(B_{3r})$ , then the estimates we can show }
	\begin{equation}\label{eq4.20}
		-\!\!\!\!\!\!\int_{B_{3r}}|Dw|^{p^+}dx\leq c\Bigg(-\!\!\!\!\!\!\int_{B_{4r}} |Dw|^{p(x)}dx+1\Bigg)
	\end{equation}
	{\em provided $r\leq \frac{R_0}{4}$}
	\begin{equation}\label{eq4.21}
		R_0\leq \frac{R}{2}\;\;\;and\;\;\;\omega(2R_0)\leq \min\Biggl\{\sqrt{\frac{n+1}{n}}-1,\log_{c_*M}2,\frac{\sigma_1}{4}\Biggr\}
	\end{equation}
	{\em where}
	\begin{equation}\label{eq4.22}
		\sigma_1:=\min\Biggl\{\frac{c_g}{16^{1/\gamma_1}},4(\gamma_1-1),1\Biggr\}
	\end{equation}
	Proof: see \cite{BCO28_72_15_13}.\vskip 1mm
	Let next note that for $x\in B_{4r}$
	\begin{equation*}
		p^+(1+\frac{\sigma_1}{4})\leq p^-(1+\frac{\sigma_1}{4})+\omega(8r)(1+\frac{\sigma_1}{4})
	\end{equation*}
	\begin{equation*}
		\hspace{8mm}\!\!\!\leq p^-(1+\frac{\sigma_1}{4})+\omega(8r)\gamma_1
	\end{equation*}
	\begin{equation}\label{eq4.23}
		\;\;\;\;\;\;\;\leq p(x)(1+\frac{\sigma_1}{4}+\omega(8r))
	\end{equation}\vskip 1mm
	Now we shall compare the weak solution $w\in W^{1,p(.)}(B_{4r})\cap W^{1,p^+}(B_{4r})$ of (\ref{eq4.2}) with a weak solution of a non-linear equation with the constant $p^+$ growth. define a function $b=b(x,\xi):B_{4r}\times\mathbb{R}^n\rightarrow \mathbb{R}$ by 
	\begin{equation}\label{eq4.24}
		b(x,\xi)=(1+|\xi|^2)^{\frac{p^+-p(x)}{2}}A(x,\xi)
	\end{equation}
	Indeed , it follows from (\ref{eqa3}) and (\ref{eq4.24}) that 
	\begin{equation*}
		|b(x,\xi)|=|(1+|\xi|^2)^\frac{p^+-p(x)}{2}A(x,\xi)|
	\end{equation*}
	\begin{equation*}
		\hspace{37mm}\leq|(1+|\xi|^2)|^\frac{p^+-p(x)}{2} \; l.|(\mu^2+|\xi|^2)|^\frac{p(x)-1}{2}
	\end{equation*}
	\begin{equation}\label{eq4.25}
		\;\;\;\leq l|(1+|\xi|^2)|^\frac{p^+-1}{2}
	\end{equation}
	By (\ref{eqh1}) and (\ref{eq4.25}), we have the following (see \cite{BOR29_64_15_5})
	\begin{equation}\label{eq4.26}
		\begin{cases}
			& (1+|\xi|^2)^\frac{1}{2}|b_\xi(x,\xi)|+|b(x,\xi)|\leq 3L'(1+|\xi|^2)^\frac{p^+-1}{2}\\
			
			& \frac{\nu'}{2}(1+|\xi|^2)^\frac{p^+-1}{2}|\eta|^2\leq \Big<b_\xi(x,\xi)\eta,\eta\Big>
		\end{cases}
	\end{equation}
	provided we choose small $R$ such that $\omega(2R)\leq \min\Big\{1,\frac{\nu'}{2L'}\Big\}$.\vskip 1mm
	Now give some known inequalities related to logarithmic functions. Let $0<\beta_1\leq\beta_2<\infty$ and $s>1$ \vskip 0.8mm
	$\bullet$ There exists $c(\beta_2,\beta_2)>0$ such that for any $\beta_1\leq \beta\leq \beta_2$ and $0<t< \mathit{e}+1$,
	\begin{equation}\label{eq4.27}
		t^\beta|\log t|\leq c(\beta_1,\beta_2)
	\end{equation}\vskip 0.1mm
	$\bullet$ There exists $c(s,\beta_1,\beta_2)>0$ such that for any $f\in L^1(\Omega)$ and $\beta\in [\beta_1,\beta_2]$ (see \cite{AM24_63_15_4})
	\begin{equation}\label{eq4.28}
		-\!\!\!\!\!\!\int_\Omega |f|\Bigg[\log\Bigg(\mathit{e}+\frac{|f|}{-\!\!\!\!\!\int_\Omega |f|dx}\Bigg)\Bigg]^\beta dx \leq c(s,\beta_1,\beta_2)\big(-\!\!\!\!\!\!\!\int_\Omega |f|^s\big)^\frac{1}{s}
	\end{equation}
	for every $f\in L\log^\beta L$. 
	
	Now the lemma as follows\vskip 1mm
	\begin{lem}\label{L4.3}
		Under the lemma and conclusion of Lemma 4.1, $\exists \; \delta=\delta(n,\nu',,L',\gamma_1,\gamma_2,\epsilon)>0$ such that if $h\in W^{1,p^+}(B_{3r})$ is the weak solution of 
		\begin{equation}\label{eq4.29}
			\begin{cases}
				& div\;b(x,Dh)=0\;\;\;in\;B_{3r}\\
				& h=w\;\;\;on\;\partial B_{3r}
			\end{cases}
		\end{equation}
		then 
		\begin{equation}\label{eq4.30}
			-\!\!\!\!\!\!\int_{B_{3r}}|Dh|^{p^+}dx\leq c\lambda
		\end{equation}
		for some constant $c=c(n,\nu',L',\gamma_1,\gamma_2)$ and 
		\begin{equation}\label{eq4.31}
			-\!\!\!\!\!\!\int_{B_{3r}} |Dw-Dh|^{p^+} dx \leq \epsilon\lambda
		\end{equation}
	\end{lem}
	{\bf Proof.} From (\ref{eq4.20}) and (\ref{eq4.4}), we have 
	\begin{equation}\label{eq4.32}
		-\!\!\!\!\!\!\int_{B_{3r}}|Dw|^{p^+} dx \leq c\lambda
	\end{equation}
	Now we test (\ref{eq4.2}) and (\ref{eq4.29}) by $h-w\in W_0^{1,p^+}(B_{3r})\subset W_0^{1,p(.)}(B_{4r})$ , to find 
	\begin{equation*}
		-\!\!\!\!\!\!\int_{B_{3r}} |Dh|^{p^+}dx\leq c\Bigg(-\!\!\!\!\!\!\!\int_{B_{3r}} |Dw|^{p^+}dx+1\Bigg)\leq c\lambda
	\end{equation*}
	Now 
	\begin{equation*}
		I_1 := -\!\!\!\!\!\!\int_{B_{3r}} \big<b(x,Dh)-b(x,Dw),D(h-w)\big>dx
	\end{equation*}
	\begin{equation*}
		\hspace{15mm}=-\!\!\!\!\!\!\int_{B_{3r}} \big<A(x,Dw)-b(x,Dw),D(h-w)\big>dx=: I_2
	\end{equation*}
	In $I_1$ , with the same mathod for the estimate (\ref{eq4.11}) , we derive that for any $k_3\in (0,1)$,
	\begin{equation}\label{eq4.33}
		-\!\!\!\!\!\!\int_{B_{3r}} |Dw-Dh|^{p^+} dx \leq c(k_3)I_1
	\end{equation}
	Now we want to estimate $I_2$, by Young's inequality there exists $k_4\in (0,1)$ such that
	\begin{equation*}
		|I_2|\leq -\!\!\!\!\!\!\int_{B_{3r}} |A(x,Dw)-b(x,Dw)||D(w-h)|dx
	\end{equation*}
	\begin{equation*}
		\hspace{38mm}\leq k_4 -\!\!\!\!\!\!\!\int_{B_{3r}}|Dw-dh|^{p^+}dx+ c(k_4) -\!\!\!\!\!\!\!\int_{B_{3r}} |A(x,Dw)-b(x,Dw)|^\frac{p^+}{p^+-1}dx
	\end{equation*}
	\begin{equation*}
		\hspace{60mm}\leq k_4  -\!\!\!\!\!\!\int_{B_{3r}} |Dw-Dh|^{p^+}dx\;+\; c(k_4) -\!\!\!\!\!\!\int_{B_{3r}} |A(x,Dw)-(1+|Dw|^2)^\frac{p^+-p(x)}{2} A(x,Dw)|^\frac{p^+}{p^+-1}dx
	\end{equation*}
	\begin{equation}\label{eq4.34}
		\hspace{60mm}\leq k_4 -\!\!\!\!\!\!\int_{B_{3r}}|Dw-Dh|^{p^+}dx\;+\; c(k_4) -\!\!\!\!\!\!\int_{B_{3r}} \Big[\underbrace{\Big\{(1+|Dw|^2)^\frac{p^+-p(x)}{2}-1\Big\}}\big|A(x,Dw)\big|\Big]^\frac{p^+}{p^+-1}dx
	\end{equation}
	Applying Mean Value Theorem (see \cite{BOR29_64_15_5}) at the underbraced term in (\ref{eq4.34}), we have for some $t_x\in [0,1]$ ,
	\begin{equation}\label{eq4.35}
		(1+|Dw|^2)^\frac{p^+-p(x)}{2}-1=\frac{p^+-p(x)}{2}(1+|Dw|^2)^{t_x\frac{p^+-p(x)}{2}}\log(1+|Dw|^2)
	\end{equation}
	Then from (\ref{eq4.35}) and (\ref{eqa3}) , I have in (\ref{eq4.34})
	\begin{equation}\label{eq4.36}
		|I_2|\leq k_4 -\!\!\!\!\!\!\!\int_{B_{3r}} |Dw-Dh|^{p^+}dx+\frac{c(k_4)}{|B_{3r}|} \int_{E}\Big[\frac{p^+-p(x)}{2}(1+|Dw|^2)^{t_x\frac{p^+-p(x)}{2}+\frac{p(x)-1}{2}}\log(1+|Dw|^2)\Big]^\frac{p^+}{p^+-1}\;dx
	\end{equation}
	As $|p(x)-p(y)|\leq \omega(|x-y|)$ and $p^+$ as defined before then we have $p^+-p(x)\leq p^+-p^-\leq \omega(8r)$, so (\ref{eq4.36}) implies 
	\begin{equation}\label{eq4.37}
		|I_2|\leq k_4 -\!\!\!\!\!\!\!\int_{B_{3r}} |Dw-Dh|^{p^+}dx+\frac{c(k_4)}{|B_{3r}|} \int_{E}\Big[\omega(8r)(1+|Dw|^2)^{t_x\frac{p^+-p(x)}{2}+\frac{p(x)-1}{2}}\log(1+|Dw|^2)\Big]^\frac{p^+}{p^+-1}\;dx
	\end{equation}
	Here $E :=\{x\in B_{3r}: \mu^2+|Dw(x)|^2>0 \}$. Since in $B_{3r}\setminus E$ i.e. when $\mu^2+|Dw(x)|^2=0$ for $x\in B_{3r}$ , then by (\ref{eqa3}) and (\ref{eq4.24}), I have $A(x,Dw)=0=b(x,Dw)$. So in $B_{3r}\setminus E$ , $I_2$ will be wiped off. So without any ambiguity we can rewrite (\ref{eq4.37}) as
	\begin{equation*}
		|I_2|\leq k_4 -\!\!\!\!\!\!\!\int_{B_{3r}} |Dw-Dh|^{p^+}dx+\frac{c(k_4)}{|B_{3r}|} \int_{B_{3r}}\Big[\omega(8r)(1+|Dw|^2)^{t_x\frac{p^+-p(x)}{2}+\frac{p(x)-1}{2}}\log(1+|Dw|^2)\Big]^\frac{p^+}{p^+-1}\;dx
	\end{equation*} 
	Now separating $B_{3r}$ into $E_1 := B_{3r}\cap \{x\in B_{3r}:|Dw(x)|\geq 1\}$ and $E_2 := B_{3r}\cap\{x\in B_{3r}:|Dw(x)|< 1\}$. By (\ref{eq4.27}) with $\frac{\gamma_1-1}{2}\leq \beta \leq \frac{\gamma_2-1}{2}$, we see that for any $x \in E_2$ ,
	\begin{equation*}
		(1+|Dw|^2)^{t_x\frac{p^+-p(x)}{2}+\frac{p(x)-1}{2}}\big|\log(1+|Dw|^2)\big|\leq c(\gamma_1,\gamma_2)
	\end{equation*}
	Thus,
	\begin{equation*}
		|I_2|\leq k_4 -\!\!\!\!\!\!\int_{B_{3r}}|Dw-Dh|^{p^+}dx\;+\;\frac{c(k_4)}{|B_{3r}|}\omega(8r)^\frac{p^+}{p^+-1}\Bigg[\int_{E_1}\Big\{(1+|Dw|^2)^{t_x\frac{p^+-p(x)}{2}+\frac{p(x)-1}{2}}\big|\log(1+|Dw|^2)\big|\Big\}^\frac{p^+}{p^+-1}\; dx \;+\; \int_{E_2} c(\gamma_1,\gamma_2)dx\Bigg]
	\end{equation*}
	\begin{equation*}
		\hspace{7mm}\leq 	k_4 -\!\!\!\!\!\!\int_{B_{3r}}|Dw-Dh|^{p^+}dx\;+\;\frac{c(k_4)}{|B_{3r}|}\omega(8r)^\frac{p^+}{p^+-1}\Bigg[\int_{E_1}\Big\{(1+|Dw|^2)^{\frac{p^+-p(x)}{2}+\frac{p(x)-1}{2}}\big|\log(1+|Dw|^2)\big|\Big\}^\frac{p^+}{p^+-1}\; dx \;+\; \int_{E_2} c(\gamma_1,\gamma_2)dx\Bigg]
	\end{equation*}
	\begin{equation*}
		\leq k_4 -\!\!\!\!\!\!\int_{B_{3r}}|Dw-Dh|^{p^+}dx\;+\;\frac{c(k_4)}{|B_{3r}|}\omega(8r)^\frac{p^+}{p^+-1}\Bigg[\int_{E_1}\Big\{(1+|Dw|^2)^{p^+}\big|\log(1+|Dw|^2)\big|\Big\}^\frac{p^+}{p^+-1}\; dx \;+\; \int_{E_2} c(\gamma_1,\gamma_2)dx\Bigg]
	\end{equation*}
	\begin{equation}\label{eq4.38}
		\leq k_4 -\!\!\!\!\!\!\int_{B_{3r}}|Dw-Dh|^{p^+}dx\;+\;c(k_4)\omega(8r)^\frac{p^+}{p^+-1}\Bigg[\underbrace{-\!\!\!\!\!\!\!\int_{B_{3r}}|Dw|^{p^+}\{\log(\mathit{e}+|Dw|^{p^+})\}^\frac{p^+}{p^+-1}\; dx }_{I_3}\;+\; 1\Bigg]
	\end{equation}
	Now we want to estimate $I_3$. Using the inequlities (see\; \cite[p. 125]{AM24_63_15_4})
	\begin{equation}\label{4.39}
		\log(\mathit{e}+ab)\leq \log(\mathit{e}+a)+\log(\mathit{e}+b) 
	\end{equation}
	and
	\begin{equation}\label{eq4.40}
		(\log(\mathit{e}+t))^\beta\leq c\alpha^{-\beta}(\mathit{e}+t)^\frac{\alpha}{4}
	\end{equation}
	for every $t>0,\alpha\in (0,1),\beta\in(\frac{\gamma_2}{\gamma_2-1},\frac{\gamma_1}{\gamma_1-1})$ and some $c=c(\gamma_1,\gamma_2)>0$
	and applying (\ref{eq4.28}) to $\beta=\frac{p^+}{p^+-1}\in [\frac{\gamma_2}{\gamma_2-1},\frac{\gamma_1}{\gamma_1-1}],s=1+\frac{\sigma_1}{4}$, we have 
	\begin{equation*}
		I_3\leq c -\!\!\!\!\!\!\int_{B_{3r}} |Dw|^{p^+}\Big[\log\Big(e+\frac{|Dw|^{p^+}}{-\!\!\!\!\!\int_{B_{3r}}|Dw|^{p^+}dx}\Big)\Big]^\frac{p^+}{p^+-1}dx\;+\; c -\!\!\!\!\!\!\!\int_{B_{3r}}|Dw|^{p^+}\Big[\log\Big({e+-\!\!\!\!\!\!\!\int_{B_{3r}}|Dw|^{p^+}dx}\Big)\Big]^\frac{p^+}{p^+-1}dx
	\end{equation*} 
	\begin{equation*}
		\;\;\;\;\leq \underbrace{c\Big(-\!\!\!\!\!\!\int_{B_{3r}}|Dw|^{p^+(1+\frac{\sigma_1}{4})}dx\Big)^\frac{1}{1+\frac{\sigma_1}{4}}}_{I_4}\;+\; c -\!\!\!\!\!\!\!\int_{B_{3r}}|Dw|^{p^+}\Big[\log\Big(e+-\!\!\!\!\!\!\int_{B_{4r}}|Dw|^{p(x)}dx+1\Big)\Big]^\frac{p^+}{p^+-1}dx\;\;\;\;\;by\;(\ref{eq4.19})
	\end{equation*}
	\begin{equation*}
		\leq I_4+ c\lambda\Big|\log\frac{1}{8r}\Big|+c\alpha^{-\frac{p^+}{p^+-1}}-\!\!\!\!\!\!\!\int_{B_{3r}}|Dw|^{p^+}\;\Big[\mathit{e}+-\!\!\!\!\!\!\int_{B_{4r}}|Dw|^{p(x)}+1\Big]^\frac{\alpha}{4}dx\;\;\;\;by\;(\ref{eq4.40})
	\end{equation*}
	\begin{equation*}
		\leq I_4 + c\lambda\Big|\log\frac{1}{8r}\Big|+c\lambda\alpha^{-\frac{p^+}{p^+-1}}\;\;\;\;\;by\;(\ref{eq4.32},\ref{eq4.4})
	\end{equation*}
	\begin{equation*}
		\leq I_4 + c\lambda\Big|\log\frac{1}{8r}\Big|+c\lambda\alpha^{-\frac{\gamma_1}{\gamma_1-1}}
	\end{equation*}
	\begin{equation}\label{eq4.41}
		\leq I_4 + c\lambda\Big|\log\frac{1}{8r}\Big|+c\lambda
	\end{equation}
	now want to estimate $I_4$ with the help of (\ref{eq4.17}),(\ref{eq4.23}) and (\ref{eq4.4})
	\begin{equation*}
		I_4\leq c\Big(-\!\!\!\!\!\!\int_{B_{3r}}(|Dw|+1)^{p(x)}dx+1\Big)^\frac{1+\frac{\sigma_1}{4}+\omega(8r)}{1+\frac{\sigma_1}{4}}
	\end{equation*}
	\begin{equation*}
		\;\;\;\;\leq c\Big(-\!\!\!\!\!\!\int_{B_{4r}}(|Dw|+1)^{p(x)}dx+1\Big)^\frac{1+\frac{\sigma_1}{4}+\omega(8r)}{1+\frac{\sigma_1}{4}}
	\end{equation*}
	\begin{equation}\label{eq4.42}
		\leq c\lambda
	\end{equation}
	Now by (\ref{eq4.41}) and (\ref{eq4.42}), (\ref{eq4.38}) implies
	\begin{equation*}
		|I_2|\leq k_4 -\!\!\!\!\!\!\!\int_{B_{3r}}|Dw-Dh|^{p^+}dx\;+\;c(k_4)\omega(8r)^\frac{p^+}{p^+-1}\Bigg[\Big|\log\frac{1}{8r}\Big|^\frac{p^+}{p^+-1}+1\Bigg]\lambda
	\end{equation*}
	\begin{equation*}
		\leq k_4 -\!\!\!\!\!\!\!\int_{B_{3r}}|Dw-Dh|^{p^+}dx\;+\;c(k_4)\Bigg[\omega(8r)\Big|\log\frac{1}{8r}\Big|+\omega(8r)\Bigg]^\frac{p^+}{p^+-1}\lambda
	\end{equation*}
	Now finally choose $R_0$ in such a way that , for every $4r\leq R_0 \leq \frac{R}{2}<\frac{1}{4}$ such that from \cite[Remark 2.3]{BCO28_72_15_13}
	\begin{equation}\label{eq4.43}
		\sup_{0<8r\leq R}\omega(8r)\leq \sup_{0<8r\leq R}\frac{\omega(8r)\log\frac{1}{8r}}{\log 2}\leq \frac{\delta}{\log 2}
	\end{equation} \vskip 1mm
	we have from (\ref{eq4.43}) and (\ref{eq2.4})
	\begin{equation} \label{eq4.44}
		|I_2|\leq k_4 -\!\!\!\!\!\!\!\int_{B_{3r}}|Dw-Dh|^{p^+}dx\;+\;c(k_4)\Bigg[\delta+\frac{\delta}{\log 2}\Bigg]^\frac{p^+}{p^+-1}\lambda
	\end{equation}
	From (\ref{eq4.44}) , (\ref{eq4.38}) implies
	\begin{equation*}
		-\!\!\!\!\!\!\int_{B_{3r}}|Dw-Dh|^{p^+}dx\leq k_4c(k_3) -\!\!\!\!\!\!\!\int_{B_{3r}}|Dw-Dh|^{p^+}dx\;+\;c(k_3)c(k_4)\Bigg[\delta+\frac{\delta}{\log 2}\Bigg]^\frac{p^+}{p^+-1}\lambda
	\end{equation*}
	\begin{equation*}
		-\!\!\!\!\!\!\int_{B_{3r}}|Dw-Dh|^{p^+}dx\leq k_4c(k_3) -\!\!\!\!\!\!\!\int_{B_{3r}}|Dw-Dh|^{p^+}dx\;+\;c(k_3)c(k_4)\Bigg[\frac{2\delta}{\log 2}\Bigg]^\frac{p^+}{p^+-1}\lambda
	\end{equation*}
	\begin{equation}\label{eq4.45}
		\hspace{32mm}\leq k_4c(k_3) -\!\!\!\!\!\!\!\int_{B_{3r}}|Dw-Dh|^{p^+}dx\;+\;c(k_3)c(k_4)\Bigg[\frac{2\delta}{\log 2}\Bigg]^\frac{p^+}{p^+-1}\lambda
	\end{equation}
	Now choose
	\begin{equation}\label{eq4.46}
		k_4=\frac{1}{2c(k_3)}\;,\;\delta<\frac{\log2}{2}\Big(\frac{\epsilon}{2c(k_3)c(k_4)}\Big)^\frac{\gamma_2-1}{\gamma_2}
	\end{equation}
	Both (\ref{eq4.45}) and (\ref{eq4.46}) implies (\ref{eq4.31}).\;\;$\Box$\vskip 2mm
	Consider
	\begin{equation}
		b_{2r}(\xi)=-\!\!\!\!\!\!\int_{B_{2r}} b(x,\xi)dx
	\end{equation}
	Then clearly $b_{B_{3r}}(\xi)$ satisfies all the conditions in (\ref{eq4.27}) with $b(x,\xi)$ replaced by $b_{B_{2r}}(\xi)$. Moreover , we have that 
	\begin{equation*}
		\frac{|b(x,\xi)-b_{B_{2r}(\xi)}|}{(\mu^2+|\xi|^2)^\frac{p_2-1}{2}}=\Bigg|\frac{b(x,\xi)}{(\mu^2+|\xi|^2)^\frac{p_2-1}{2}}-\frac{-\!\!\!\!\!\int_{B_{2r}}b(x,\xi)dx}{(\mu^2+|\xi|^2)^\frac{p_2-1}{2}}\Bigg|=\Bigg|\frac{A(x,\xi)}{(\mu^2+|\xi|^2)^\frac{p(x)-1}{2}}--\!\!\!\!\!\!\int_{B_{2r}}{\frac{A(x,\xi)dx}{(\mu^2+|\xi|^2)^\frac{p(x)-1}{2}}}dx\Bigg|
	\end{equation*}
	Since $a(x,\xi)$ satisfies (\ref{eq2.7}), so is $b(x,\xi)$. Before going to discuss a lemma  Now we would like to recall a result related to provide Lipschitz regularity of weak solutions for non-linear elliptical equations under constant exponent $p$ and the nonlinearity independent of $x$ variables. For finte $p(>1)$, let $\boldsymbol{b}=\boldsymbol{b}(\xi) : \mathbb{R}^n\rightarrow \mathbb{R}^n$ be a vector values function satisfying 
	\begin{equation}\label{eq4.48}
		\begin{cases}
			& (1+|\xi|^2)^\frac{1}{2}|\boldsymbol{b}_\xi(\xi)|+|\boldsymbol{b}(\xi)|\leq L(1+|\xi|^2)^\frac{p-1}{2}\\
			
			& \nu(1+|\xi|^2)^\frac{p-1}{2}|\eta|^2\leq \Big<\boldsymbol{b}_\xi(\xi)\eta,\eta\Big>
		\end{cases}
	\end{equation}
	ffor all $\xi,\eta\in\mathbb{R}^n$ and for positive constants $\nu\leq L<+\infty$. Then we have the following result
	\begin{lem}\label{L4.4}(see \cite[Lemma 3.2(i)]{BCO28_72_15_13},\cite{L54_52p_21},\cite[Theorem 2.1]{L55_52p_22})
		If $v\in W^{1,p}(B_{2r})$ is a weak solution of 
		\begin{equation*}
			div\;\boldsymbol{b}(Dv)=0\;\;\;in\;B_{2r}
		\end{equation*}
		then $Dv\in L^\infty(B_r)$with the estimate
		\begin{equation*}
			||Dv||^p_{L^\infty(B_r)}\leq c\Bigg( -\!\!\!\!\!\!\int_{B_{2r}}|Dv|^pdx+1\Bigg)
		\end{equation*}
		for constant $c=c(n,\nu',L',p)$.
	\end{lem}
	Now discuss a lemma
	\begin{lem}\label{L4.5}
		Under the assumptions and conclusion of previous Lemmas \ref{L4.1},\ref{L4.2} and \ref{L4.3} , there exists $\delta=\delta(n,\nu',L',\gamma_1,\gamma_2,\epsilon)$ such that if
		\begin{equation*}
			-\!\!\!\!\!\!\int_{B_{3r}}|Dh|^{p^+}dx\leq \lambda
		\end{equation*}
		and $v\in W^{1,p^+}(B_2r)$ is the unique solution of 
		\begin{equation}\label{eq4.49}
			\begin{cases}
				& div\;b_{B_{2r}}(Dv)=0\;\;\;in\;B_{2r}\\
				& v=h\;\;\;on\;\partial B_{2r}
			\end{cases}
		\end{equation}
		then 
		\begin{equation*}
			-\!\!\!\!\!\!\int_{B_{2r}}|Dv|^{p^+}dx\leq c\lambda
		\end{equation*}
		for some $c=c(n,\nu',L',\gamma_1,\gamma_2)$ and 
		\begin{equation*}
			-\!\!\!\!\!\!\int_{B_{2r}}|Dh-Dv|^{p^+}dx\leq \epsilon\lambda
		\end{equation*}
	\end{lem}   
	{\bf Proof.} see \cite[Lemma 4.4]{BOR29_64_15_5}.\vskip 1mm
	So from the previous lemmas (\ref{L4.1})-(\ref{L4.5}), we 
	have finally the following comparison results
	\begin{lem}\label{L4.6}
		Let $\epsilon\in (0,1)$ and $\lambda\geq 1$. There exists a constants $\delta=\delta(n,\nu',L',\gamma_1,\gamma_2,\epsilon)$, $p(.)$ and $A$ satisfies {\em log-Holder} continuity result (\ref{eq2.4}) and (\ref{eq2.6}) respectively, and  if
		\begin{equation*}
			-\!\!\!\!\!\!\int_{B_{4r}}(|Du|^{p(x)}+|u|{^p(x)})dx\leq \lambda\;\;\;and\;\;\;-\!\!\!\!\!\!\int_{B_{4r}}\Psi^{p(x)}dx\leq \delta\lambda
		\end{equation*} 
		for some $B_{4r}\subset\Omega$ with $r\leq R_0/4$, then I have $w\in W^{1,p^+}(B_r)$ and $v\in W^{1,\infty}(B_r)$ such that for any $p>1$
		\begin{equation*}
			-\!\!\!\!\!\!\int_{B_{r}}|Du-Dw|^{p(x)}dx\leq\epsilon\lambda\;,\;-\!\!\!\!\!\!\int_{B_{r}}|Dw-Dv|^{p^+}dx\leq \epsilon\lambda\;\;and\;\;||Dv||^p_{L^\infty(B_r)}\leq c_1\lambda
		\end{equation*}
	\end{lem}
	\subsection{Proof of Theorem \ref{T2.1}}
	In this section , we shall prove Calder\'on-Zygmund type estimates.These following calculations in this subsection are mainly based on the ideas from\cite{G20_15_1_24} and \cite{BCO28_72_15_13}. For a given $a(x)^\frac{p(x)}{p(x)-1}$,$|D\psi(x)|^{p(x)}\in L^q_{loc}(\Omega)$. we also assume that $p(.)$ and $A$ satisfies {\em log-H$\ddot{o}$lder} continuity result (\ref{eq2.4}) and (\ref{eq2.6}) respectively , $\delta$ being determined while $R\leq \frac{1}{2}$ being an arbitrary given small number. Further assume that $R_0>0$ satisfies 
	\begin{equation*}
		R_0\leq\frac{R}{2}\;\;\;\;\;and\;\;\;\;\; \omega(2R_0)\leq \min\Bigg\{\sqrt{\frac{n+1}{n}}-1,\log_{c_*M}2,\frac{\sigma_1}{4}\Bigg\}
	\end{equation*}
	where $\sigma_1$ and $c_*M$ are defined as (\ref{eq4.22}) and (\ref{eq4.19}) respectively. By our assumption for any relatively compact set $K\subset\Omega$
	\begin{equation*}
		\int_K(|D\psi|^{p(x)})^q<\infty\implies\int_K(|D\psi|^q)^{p(x)}<\infty
	\end{equation*}
Since $p(.)$ is bounded on $\Omega$ so $|\Omega_\infty|=0$ then by argument in subsection \ref{s3.1} we can say $|D\psi|^q\in L^{p(.)}(\Omega)$. Then by Proposition \ref{p3.1}, we can say $|D\psi|^q\in L^1_{loc}(\Omega)$, so $|D\psi|\in L^q_{loc}(\Omega)$. Now observe that the inclusion
\begin{equation*}
	L^{p(.)q}_{loc}(\Omega)\subset L^q_{loc}(\Omega)\subset L^{\gamma_2\tilde{q},\tilde{\lambda}}_{loc}(\Omega)
\end{equation*} 
with $\tilde{q}>1$,$n-\gamma_1<\tilde{\lambda}<n$ and only if $q>n$. \vskip 1mm
	Then, for the weak solution $u\in W^{1,p(.)(\Omega)}$ of (\ref{eq1.3}) and for the function $\Psi= a(x)^\frac{1}{p(x)-1}+|D\psi(x)|+1$, we define
	\begin{equation*}
		E_R(|Du|^{p(.)},\lambda)=\{x\in B_R : |Du(x)|^{p(x)}+|u(x)|^{p(x)}>\lambda\}\;\;,\;\;E_R(\Psi^{p(.)},\lambda)=\{x\in B_R : \Psi(x)^{p(x)}>\lambda\}
	\end{equation*}
	for any $\lambda>0$. Let set
	\begin{equation}\label{eq4.50}
		\lambda_0 = -\!\!\!\!\!\!\int_{B_{R}}\Big(|Du|^{p(x)}+|u|^{p(x)}+\frac{1}{\delta}\Psi^{p(x)}\Big)dx
	\end{equation}
	For any $y\in B_R$ ,  define a continuous map $G(y,.): (0,R_0]\rightarrow [0,\infty)$ by 
	\begin{equation*}
		G(y,r)=-\!\!\!\!\!\!\int_{B_{r}(y)}\Big(|Du|^{p(x)}+|u|^{p(x)}+\frac{1}{\delta}\Psi^{p(x)}\Big)dx
	\end{equation*}
	So, for almost every $y\in E_R(|Du|^{p(.)},\lambda)$, it follows from the Lebesgue differentiation theorem (see \cite[Appendix E]{E53_127}) that 
	\begin{equation}\label{eq4.51}
		\lim_{r\rightarrow 0}G(y,r)=|Du|^{p(x)}+|u|^{p(x)}+\frac{1}{\delta}\Psi^{p(x)}>\lambda
	\end{equation} 
	On the other hand , observe that for any $r\in [\frac{R}{L^{1/n}},R_0]$, and by (\ref{eq4.50})
	\begin{equation}\label{eq4.52}
		G(y,r)\leq \frac{|B_R|}{|B_r(y)|}-\!\!\!\!\!\!\int_{B_{R}}\Big(|Du|^{p(x)}+|u|^{p(x)}+\frac{1}{\delta}\Psi^{p(x)}\Big)dx\leq \frac{R^n}{r^n}\lambda_0\leq L\lambda_0<\lambda
	\end{equation}
	for $\lambda$ and $L$ selected sufficiently large such that
	\begin{equation}\label{eq4.53}
		\lambda>\lambda_1 := L\lambda_0\;\;\;and \;\;\;L\geq \Bigg(\frac{R}{R_0}\Bigg)^n
	\end{equation}
	Now fix any $\lambda>\lambda_1$. Since $G(,.)$ is a continuous map, estimates (\ref{eq4.51}) and (\ref{eq4.52}) imply that for almost every $y\in E_R(|Du|^{p(.)},\lambda)$ there exists a number $r_y\in (0,\frac{R}{L^{1/n}})$ such that
	\begin{equation*}
		G(y,r_y)=\lambda\;\;\;and\;\;\;G(y,r)<\lambda\;\;\;\forall r\in(r_y,R_0]
	\end{equation*}
	In view of the {\em Vitali covering lemma} , there is a countable family of mutually disjoints balls $\{B_{r_i}(y^i)\}^\infty_{i=1}$ with $y^i\in E_R(|Du|^{p(.)},\lambda)$ and $r_i\in (0,\frac{R}{L^{1/n}})$ such that 
	\begin{equation}\label{eq4.54}
		E_R(|Du|^{p(.)},\lambda)\subset \bigcup_{i\geq 1} B_{{5r}_i}(y^i) \cup (a\;negligible\;set)
	\end{equation}
	\begin{equation}\label{eq4.55}
		-\!\!\!\!\!\!\int_{B_{{r}_i}(y^i)}\Big(|Du|^{p(x)}+|u|^{p(x)}+\frac{1}{\delta}\Psi^{p(x)}\Big)dx=\lambda
	\end{equation}
	\begin{equation}\label{eq4.56}
		-\!\!\!\!\!\!\int_{B_{{r}_i}(y^i)}\Big(|Du|^{p(x)}+|u|^{p(x)}+\frac{1}{\delta}\Psi^{p(x)}\Big)dx<\lambda\;\;\;\;for \;\;each\;r\in(r_i,R_0]
	\end{equation}
	Now let fix $i$. From (\ref{eq4.55}) have
	\begin{equation}\label{eq4.57}
		\lambda=\frac{1}{|B_{r_i}(y^i)|}\int_{B_{{r}_i}(y^i)}\Big(|Du|^{p(x)}+|u|^{p(x)}+\frac{1}{\delta}\Psi^{p(x)}\Big)dx
	\end{equation}
	and then 
	\begin{equation}\label{eq4.58}
		|B_{r_i}(y^i)|=\frac{1}{\lambda}\int_{B_{{r}_i}(y^i)}\Big(|Du|^{p(x)}+|u|^{p(x)}\Big)dx+\frac{1}{\lambda}\int_{B_{{r}_i}(y^i)}\frac{\Psi^{p(x)}}{\delta}dx=\frac{1}{\lambda}\boldsymbol{A_1}+\frac{1}{\delta\lambda}\boldsymbol{A_2}
	\end{equation}
	Now want to estimate integral $\boldsymbol{A_1}$ as follows
	\begin{equation*}
		\boldsymbol{A_1}=\int_{{B_{{r}_i}(y^i)}\cap\{|Du|^{p(x)}+|u|^{p(x)}\leq \frac{\lambda}{4}\}}(|Du|^{p(x)}+|u|^{p(x)})dx\;+\;\int_{{B_{{r}_i}(y^i)}\cap\{|Du|^{p(x)}+|u|^{p(x)}> \frac{\lambda}{4}\}}(|Du|^{p(x)}+|u|^{p(x)})dx
	\end{equation*}
	\begin{equation}\label{eq4.59}
		\leq \frac{\lambda}{4}|B_{r_i}(y^i)|+\int_{{B_{{r}_i}(y^i)}\cap\{|Du|^{p(x)}+|u|^{p(x)}> \frac{\lambda}{4}\}}(|Du|^{p(x)}+|u|^{p(x)})dx
	\end{equation}
	In a similar way, we shall get the following estimate for the integral $\boldsymbol{A_2}$
	\begin{equation*}
		\boldsymbol{A_2}=\int_{{B_{{r}_i}(y^i)}\cap\{\Psi^{p(x)}\leq \frac{\delta\lambda}{4}\}}\Psi^{p(x)}dx\;+\;\int_{{B_{{r}_i}(y^i)}\cap\{\Psi^{p(x)}>\frac{\delta\lambda}{4}\}}\Psi^{p(x)}dx
	\end{equation*}
	\begin{equation}\label{eq4.60}
		\leq\frac{\delta\lambda}{4}|B_{r_i}(y_i)|+\int_{{B_{{r}_i}(y^i)}\cap\{\Psi^{p(x)}>\frac{\delta\lambda}{4}\}}\Psi^{p(x)}dx
	\end{equation}
	Inserting estimates (\ref{eq4.59}) and (\ref{eq4.60}) in (\ref{eq4.58}) we obtain
	\begin{equation*}
		\frac{|B_{r_i}(y^i)|}{2}\leq\frac{1}{\lambda}\Bigg(\int_{{B_{{r}_i}(y^i)}\cap\{|Du|^{p(x)}+|u|^{p(x)}> \frac{\lambda}{4}\}}(|Du|^{p(x)}+|u|^{p(x)})dx+\frac{1}{\delta}\int_{{B_{{r}_i}(y^i)}\cap\{\Psi^{p(x)}>\frac{\delta\lambda}{4}\}}\Psi^{p(x)}\Bigg)
	\end{equation*}
	Thus we get
	\begin{equation}\label{eq4.61}
		|B_{r_i}(y^i)|\leq\frac{2}{\lambda}\Bigg(\int_{{B_{{r}_i}(y^i)}\cap E_R(|Du|^{p(x)},\frac{\lambda}{4}) }(|Du|^{p(x)}+|u|^{p(x)})dx+\frac{1}{\delta}\int_{{B_{{r}_i}(y^i)}\cap E_R(\Psi^{p(x)},\frac{\delta\lambda}{4})}\Psi^{p(x)}\Bigg)
	\end{equation}
	Let $\epsilon\in (0,1)$. Choosing $L\geq\max\{20^n,
	(R/R_0)^n\}$. Then $r_i<\frac{R}{L^{1/n}}\leq \frac{R}{20}$ implies $20r_i\leq R$ . Thus we have $B_{20r_i}(y^i)\subset\Omega$ . In this case , from (\ref{eq4.56}),we get
	\begin{equation}\label{eq4.62}
		-\!\!\!\!\!\!\int_{B_{20r_i}(y^i)}(|Du|^{p(x)}+|u|^{p(x)})dx \leq \lambda\;\;\;and\;\;\; -\!\!\!\!\!\!\!\int_{B_{20r_i}(y^i)}\Psi^{p(x)}dx\leq \delta\lambda
	\end{equation}
	according to last lemma \ref{L4.6} , I can find sufficiently small $\delta=\delta(n,\nu,L,l',\gamma_1,\gamma_2,\epsilon),w_i\in W^{1,p(.)}(B_{5r_i}(y^i))$ and $v_i\in W^{1,\infty}(B_{5r_i}(y^i))$ satisfying
	\begin{equation}\label{eq4.63}
		-\!\!\!\!\!\!\int_{B_{5r_i}(y^i)}|Du-Dw_i|^{p(x)}dx\leq \epsilon\lambda\;,\;-\!\!\!\!\!\!\int_{B_{5r_i}(y^i)}|Dw_i-Dv_i|^{p_2^i}dx\leq \epsilon\lambda
	\end{equation} 
	and 
	\begin{equation}\label{eq4.64}
		||Dv_i||^{p_2^i}_{L^\infty(B_{5r_i}(y^i))}\leq c_1\lambda
	\end{equation}
	for some $c_1=c_1(n,\nu,L,l',\gamma_1,\gamma_2)$ , independent of $i$ and $\lambda$ and $p_2^i=\sup_{x\in B_{20r_i}(y^i)}p(x)$.\vskip 1mm
	Now we can prove Calder\'on-Zygmund related theorem \\[1.5mm]
	{\bf Proof.} Let $c_2=2.4^{\gamma_2-1}(c_1+2)>1$,where $c_1$ is given in (\ref{eq4.64}). For each $\lambda\geq \lambda_1$, we recall the covering $\{B_{r_i}(y^i)\}$ satisfying (\ref{eq4.54})-(\ref{eq4.56}), I define
	\begin{equation*}
		\tilde{E}_R(|Du|^{p(.)},\lambda)=\{x\in B_R:|Du(x)|^{p(x)}>\lambda\}
	\end{equation*}
	Since I have $\tilde{E}_R(|Du|^{p(.)},\lambda)\subset E_R(|Du|^{p(.)},c_2\lambda)\subset E_R(|Du|^{p(.)},\lambda)$ we can conclude that 
	\begin{equation}\label{eq4.65}
		\int_{\tilde{E}_R(|Du|^{p(.)},c_2\lambda)}|Du|^{p(x)}\leq \sum_{i=1}^{\infty}\Bigg(\int_{\tilde{E}_R(|Du|^{p(.)},c_2\lambda)\cap B_{5r_i}(y^i)}|Du|^{p(x)}dx\Bigg)
	\end{equation}
	For $B_{20r_i(y^i)}\subset\Omega$, use the definition of $c_2$ and (\ref{eq4.64}) , using elementary computation $(a+b)^p\leq 2^{p-1}(a^p+b^p)$ for $p>1$ to find that for a.e. $x\in \tilde{E}_R(|Du|^{p(.)},c_2\lambda)\cap B_{5r_i}(y^i)$, 
	\begin{equation*}
		|Du|^{p(x)}=|Du-Dw_i+Dw_i-Dv_i+Dv_i|^{p(x)}
	\end{equation*}
	\begin{equation*}
		\hspace{39mm}\leq 4^{\gamma_2-1}(|Du-Dw_i|^{p(x)}+|Dw_i-Dv_i|^{p(x)}+|Dv_i|^{p(x)})
	\end{equation*}
	\begin{equation*}
		\hspace{44mm}\leq 4^{\gamma_2-1}(|Du-Dw_i|^{p(x)}+|Dw_i-Dv_i|^{p_2^i}+|Dv_i|^{p^i_2}(x)+2)
	\end{equation*}
	\begin{equation*}
		\hspace{37mm}\leq 4^{\gamma_2-1}(|Du-Dw_i|^{p(x)}+|Dw_i-Dv_i|^{p_2^i}+(c_1+2)\lambda)
	\end{equation*}
	\begin{equation*}
		\hspace{37mm}\leq 4^{\gamma_2-1}(|Du-Dw_i|^{p(x)}+|Dw_i-Dv_i|^{p_2^i})+\frac{|Du|^{p(x)}}{2}
	\end{equation*}
	Thus
	\begin{equation*}
		\;\;\;\;\;\;\;\;|Du|^{p(x)}\leq 2.4^{\gamma_2-1}(|Du-Dw_i|^{p(x)}+|Dw_i-Dv_i|^{p_2^i})
	\end{equation*}
	Then , from (\ref{eq4.63}) we have
	\begin{equation*}
		\int_{\tilde{E}_R(|Du|^{p(.)},c_2\lambda)\cap B_{5r_i}(y^i)}|Du|^{p(x)}dx\leq 2.4^{\gamma_2-1}\Bigg(\int_{B_{5r_i}(y^i)}|Du-Dw_i|^{p(x)}+\int_{B_{5r_i}(y^i)}|Dw_i-Dv_i|^{p_2^i}dx\Bigg)
	\end{equation*}
	\begin{equation}\label{eq4.66}
		\leq 4.4^{\gamma_2-1}5^n|B_{r_i}(y^i)|\epsilon\lambda
	\end{equation}
	Therefore,combining (\ref{eq4.66}) and (\ref{eq4.61}) , we have for each $i$ ,
	\begin{equation}\label{eq4.67}
		\int_{\tilde{E}_R(|Du|^{p(.)},c_2\lambda)\cap B_{5r_i}(y^i)}|Du|^{p(x)}dx\leq C_3\epsilon \Bigg(\int_{{B_{{r}_i}(y^i)}\cap E_R(|Du|^{p(x)},\frac{\lambda}{4}) }(|Du|^{p(x)}+|u|^{p(x)})dx+\frac{1}{\delta}\int_{{B_{{r}_i}(y^i)}\cap E_R(\Psi^{p(x)},\frac{\delta\lambda}{4})}\Psi^{p(x)}\Bigg)
	\end{equation}
	where $c_3=8.4^{\gamma_2-1}5^n$. Since $\{B_{r_i}(y^i)\}$ are mutually disjoints, summing on $i$, it follows from (\ref{eq4.65}) and (\ref{eq4.67}) that
	\begin{equation}\label{eq4.68}
		\int_{\tilde{E}_R(|Du|^{p(.)},c_2\lambda)}|Du|^{p(x)}dx\leq c_3\epsilon \Bigg(\int_ {E_R(|Du|^{p(x)},\frac{\lambda}{4}) }(|Du|^{p(x)}+|u|^{p(x)})dx+\frac{1}{\delta}\int_ {E_R(\Psi^{p(x)},\frac{\delta\lambda}{4})}\Psi^{p(x)}dx\Bigg)
	\end{equation}
	for any $\lambda\geq \lambda_1$.\vskip 1mm
	Now introduce the truncated function/cutting operator (see for references \cite[p. 28]{BL52},\cite[p. 16]{O_55_129})
	\begin{equation*}
		T_k(\sigma)=\min\{k,\sigma\}
	\end{equation*}
	for every $k,\sigma\in\mathbb{R}$ and define for every $\tilde{\lambda}>0$ 
	\begin{equation*}
		F_k(|Du|^{p(.)},\tilde{\lambda})=\{x\in B_R : T_k(|Du(x)|^{p(x)}+|u(x)|^{p(x)})>\tilde{\lambda}\}
	\end{equation*}
	and
	\begin{equation*}
		\tilde{F}_k(|Du|^{p(.)},\tilde{\lambda})= \{x\in B_R : T_k(|Du(x)|^{p(x)})>\tilde{\lambda}\}
	\end{equation*}
	Notice that (\ref{eq4.68}) implies
	\begin{equation}\label{eq4.69}
		\int_{\tilde{F}_k(|Du|^{p(.)},c_2\lambda)}|Du|^{p(x)}dx\leq c_3\epsilon\Bigg(\int_ {F_k(|Du|^{p(.)},\frac{\lambda}{4})}(|Du|^{p(x)}+|u|^{p(x)})dx+\frac{1}{\delta}\int_ {F_k(\Psi^{p(.)},\frac{\delta\lambda}{4})}\Psi^{p(x)}dx\Bigg)
	\end{equation}
	Multiplying (\ref{eq4.69}) with $\lambda^{q-2}$ and integrating with respect to $\lambda$ for $\lambda\in(\lambda_1,+\infty)$, we get
	\begin{equation*}
		\int_{\lambda_1}^{+\infty}\lambda^{q-2}\Bigg(\int_{\tilde{F}_k(|Du|^{p(.)},c_2\lambda)}|Du|^{p(x)}dx\Bigg)d\lambda\leq c_3\epsilon\int_{\lambda_1}^{+\infty}\lambda^{q-2}\Bigg(\int_ {F_k(|Du|^{p(.)},\frac{\lambda}{4})}(|Du|^{p(x)}+|u|^{p(x)})dx\Bigg)d\lambda
	\end{equation*}
	\begin{equation}\label{eq4.70}
		\hspace{60mm}+\;c_3\epsilon\int_{\lambda_1}^{+\infty}\lambda^{q-2}\Bigg(\int_{F_k(\Psi^{p(.)},\frac{\delta\lambda}{4})}\frac{\Psi^{p(x)}}{\delta}\;dx\Bigg)d\lambda
	\end{equation}
	Now we estimate the left hand side of (\ref{eq4.70}). By using Fubini's Theorem, we obtain
	\begin{equation*}
		\int_{\lambda_1}^{+\infty}\lambda^{q-2}\Bigg(\int_{\tilde{F}_k(|Du|^{p(.)},c_2\lambda)}|Du|^{p(x)}dx\Bigg)d\lambda=\int_{\tilde{F}_k(|Du|^{p(.)},c_2\lambda)}|Du|^{p(x)}\Bigg(\int_{\lambda_1}^{\frac{T_k(|Du|^{p(.)})}{c_2}}\lambda^{q-2}\Bigg)dx
	\end{equation*}
	\begin{equation*}
		=\frac{1}{q-1}\int_{\tilde{F}_k(|Du|^{p(.)},c_2\lambda)}|Du|^{p(x)}\Biggl\{\frac{[T_k(|Du|^{p(x)})]^{q-1}}{c_2^{q-1}}-\lambda_1^{q-1}\Biggr\}dx
	\end{equation*}
	\begin{eqnarray*}
		=\frac{1}{q-1}\Biggl\{\int_{B_R}\frac{|Du|^{p(x)}[T_k(|Du|^{p(x)})]^{q-1}}{c_2^{q-1}}dx-\int_{T_k(|Du|^{p(x)})\leq c_2\lambda}\frac{|Du|^{p(x)}[T_k(|Du|^{p(x)})]^{q-1}}{c_2^{q-1}}-\int_{B_{R}}\lambda_1^{q-1}|Du|^{p(x)}dx\;+\\
		\int_{T_k(|Du|^{p(x)})\leq c_2\lambda}\lambda_1^{q-1}|Du|^{p(x)}\Biggr\}
	\end{eqnarray*}
	\begin{eqnarray*}
		=\frac{1}{q-1}\Biggl\{\int_{B_R}\frac{|Du|^{p(x)}[T_k(|Du|^{p(x)})]^{q-1}}{c_2^{q-1}}dx-\int_{T_k(|Du|^{p(x)})\leq c_2\lambda}\frac{|Du|^{p(x)}c_2^{q-1}\lambda_1^{q-1}}{c_2^{q-1}}-\int_{B_{R}}\lambda_1^{q-1}|Du|^{p(x)}dx\;+\\
		\int_{T_k(|Du|^{p(x)})\leq c_2\lambda}\lambda_1^{q-1}|Du|^{p(x)}\Biggr\}
	\end{eqnarray*}
	\begin{equation}\label{eq4.71}
		=\frac{1}{q-1}\Biggl\{\int_{B_R}\frac{|Du|^{p(x)}[T_k(|Du|^{p(x)})]^{q-1}}{c_2^{q-1}}dx-\int_{B_R}\lambda_1^{q-1}|Du|^{p(x)}dx\Biggr\}
	\end{equation}
	Similarly, we can estimate the integrals in the right hand side of (\ref{eq4.70}). Thus we can get
	\begin{equation}\label{eq4.72}
		c_3\epsilon \int_{\lambda_1}^{+\infty}\lambda^{q-2}\Bigg(\int_ {F_k(|Du|^{p(.)},\frac{\lambda}{4})}(|Du|^{p(x)}+|u|^{p(x)})dx\Bigg)d\lambda \leq c_3\epsilon\frac{4^{q-1}}{q-1}\int_{B_R}(|Du|^{p(x)}+|u|^{p(x)})[T_k(|Du|^{p(x)}+|u|^{p(x)})]^{q-1}
	\end{equation}
	\begin{equation}\label{eq4.73}
		c_3\epsilon \int_{\lambda_1}^{+\infty}\lambda^{q-2}\Bigg(\int_{F_k(\Psi^{p(.)},\frac{\delta\lambda}{4})}\frac{\Psi^{p(x)}}{\delta}\;dx\Bigg)d\lambda\leq c_3\epsilon\frac{4^{q-1}}{\delta^{q-1}}\int_{B_R}\Psi^{p(x)q}
	\end{equation}
	Consider $s,t\geq 0$. If $s\leq t$, then $T_k(s)\leq T_k(t)$. Moreover, I have 
	\begin{equation*}
		(s+t)T_k(s+t)\leq
		\begin{cases}
			& 2sT_k(2s),\;\;if\;t\leq s\\
			& 2tT_k(2t),\;\;if\;s\leq t
		\end{cases}
	\end{equation*} 
	that implies
	\begin{equation}\label{eq4.74}
		(s+t)T_k(s+t)\leq 4(sT_k(s)+tT_k(t))
	\end{equation}
	So, by virtue of (\ref{eq4.71})-(\ref{eq4.74}), we obtain
	\begin{eqnarray*}
		\frac{1}{q-1}\int_{B_R}\frac{|Du|^{p(x)}[T_k(|Du|^{p(x)})]^{q-1}}{c_2^{q-1}}dx\leq \frac{1}{q-1}\int_{B_R}\lambda_1^{q-1}|Du|^{p(x)}dx+c_3\epsilon\frac{4^{q-1}}{q-1}\int_{B_R}|Du|^{p(x)}[T_k(|Du|^{p(x)})]^{q-1}\\
		+\;c_3\epsilon\frac{4^{q-1}}{q-1}\int_{B_R}|u|^{p(x)}[T_k(|u|^{p(x)})]^{q-1}+ c_3\epsilon\frac{4^{q-1}}{\delta^{q-1}}\int_{B_R}\Psi^{p(x)q}dx
	\end{eqnarray*}
	Choosing $\epsilon$ sufficiently small, we get
	\begin{equation*}
		\int_{B_R}|Du|^{p(x)}[T_k(|Du|^{p(x)})]^{q-1}dx\leq c\int_{B_R}|Du|^{p(x)}dx\;+\;c\int_{B_R}|u|^{p(x)}[T_k(|u|^{p(x)})]^{q-1}\;+\; c \int_{B_R}\Psi^{p(x)q}dx
	\end{equation*}
	By letting $k\rightarrow\infty$, we get
	\begin{equation*}
		\int_{B_R}|Du|^{p(x)q}dx\leq c\int_{B_R}\Big(|Du|^{p(x)}+|u|^{p(x)q}+\Psi^{p(x)q}\Big)dx
	\end{equation*}
	This completes the proof.\;\;\;$\Box$\vskip 3mm
	\section{Higher Differentiability}
	Before going to prove the theorem \ref{T2.2} , let now discuss some important notations. Let take $\mathfrak{B}\Subset\Omega$, a ball of radius $R$. we denote $\mathcal{Q}_{inn}(\mathfrak{B})$ and $\mathcal{Q}_{out}(\mathfrak{B})$ as the largest and smallest cubes, cocentric to $\mathfrak{B}$, contained in and containinig $\mathfrak{B}$ respectively. Sides of $\mathcal{Q}_{inn}(\mathfrak{B})$ and $\mathcal{Q}_{out}(\mathfrak{B})$ are parallel to coordinate axes. Then clearly $|\mathfrak{B}|\approx|\mathcal{Q}_{inn}|\approx|\mathcal{Q}_{out}|$. Let denote enlarged ball as $\hat{\mathfrak{B}}=8\mathfrak{B}$. Then we have following set inclusions
	\begin{equation*}
		\mathcal{Q}_{inn}(\mathfrak{B})\subset\mathfrak{B}\subset 2\mathfrak{B}\Subset 4\mathfrak{B}\Subset\mathcal{Q}_{inn}(\hat{\mathfrak{B}})\subset\mathfrak{B}\subset \mathcal{Q}_{out}(\hat{\mathfrak{B}})
	\end{equation*}
	provided we always assume that we have choosen $\mathfrak{B}$ in such a way that $\mathcal{Q}_{out}(\hat{\mathfrak{B}})\Subset\Omega$.\vskip 1mm
	Now select $R_0$ as in (\ref{eq4.21}) such that $\omega(8r)\leq \omega(2R_0)\leq \frac{\sigma_1}{4}$. Consider the ball $B_{4r}\Subset B_{R_0}\Subset B_{R/2}$ , I choose $p^+$ and $p^-$ as before , then by (\ref{eq4.21}) and (\ref{eq4.23})
	\begin{equation}\label{eq5.1}
		p^+(1+\frac{\sigma_1}{4})\leq p(x)(1+\frac{\sigma_1}{2})
	\end{equation}
	Now we will state a lemma frequently will be used in the sequel
	\begin{lem}\label{L5.1}
		Let $u\in\mathcal{K}_\psi(\Omega)$ be a solution of obstacle variational problem satisfying (\ref{eq1.3}), $p^+,R,\delta$ defined as above and the conditions assumed in theorem \ref{T2.1}, then for any $r_0<R$ and $q>\frac{3}{2}$, I have 
		\begin{equation*}
			Du\in L^{p^+}(B_{r_0})
		\end{equation*}
	\end{lem}
	{\bf Proof.} From estimates in theorem \ref{T2.1},(\ref{eq4.22})  and (\ref{eq5.1})
	\begin{equation*}
		\int_{B_{r_0}}|Du|^{p^+}dx\leq \int_{B_{r_0}}(1+|Du|)^{p^+}dx\leq \int_{B_{r_0}}(1+|Du|)^{p(x)(1+\frac{\sigma_1}{2})}<+\infty
	\end{equation*}
	After this preliminery discussion , let start to prove theorem \ref{T2.2}\vskip 1mm
	After this preliminary observations , let start proving higher differentiability result.\\
	Let consider $\phi := u+tv$ for a appropriate $v\in W^{1,p(x)}_0(\Omega)$ such that 
	\begin{equation}\label{eq5.2}
		u-\psi+tv\geq 0\;\;\;for\;t\in[0,1)
	\end{equation}
	It is obvious that such function $\phi$ belongs to the obstacle class $\mathcal{K}_\psi(\Omega)$, because $\phi=u+tv\geq \psi$.\vskip 0.5mm
	Consider a cut off function $\eta\in C^\infty_0(B_r)$, $\eta\equiv 1$ on $B_{\frac{r}{2}}$ such that $|\nabla\eta|\leq \frac{c}{r}$. Then , for $|h|<\frac{r}{4}$, we consider
	\begin{equation}\label{eq5.3}
		v_1(x)=\eta^2(x)[(u-\psi)(x+h)-(u-\psi)(x)]
	\end{equation}
	It is easy to check that $v_1\in W^{1,p(x)}_0(\Omega)$ from the regularity of $u$.and $\psi$. Moreover $v_1$ fulfills (\ref{eq5.2}).\\
	Because of $u\in\mathcal{K}_\psi(\Omega)$ for a.e. $x\in\Omega$ and for any $t\in [0,1)$ we have 
	\begin{equation*}
		u(x)-\psi(x)+tv_1(x)
	\end{equation*}
	\begin{equation*}
		=u(x)-\psi(x)+t\eta^2(x)[(u-\psi)(x+h)-(u-\psi)(x)]
	\end{equation*}
	\begin{equation*}
		=t\eta^2(x)[(u-\psi)(x+h)+(1-t\eta^2(x))(u-\psi)(x)]\geq 0
	\end{equation*}
	\vskip 0.8mm
	By using $\phi=u+tv_1$ in (\ref{eq1.3}) as an admissible test function, we obtain
	\begin{equation}\label{eq5.4}
		\int_\Omega\big<\mathcal{A}(x,u(x),Du(x)),D(\eta^2(x)\tau_h(u-\psi))\big>dx\geq \int_\Omega\mathcal{B}(x,u(x),Du(x)).[\eta^2(x)\tau_h(u-\psi)]dx 
	\end{equation}
	On the other hand , If we define 
	\begin{equation}\label{eq5.5}
		v_2(x)=\eta^2(x-h)[(u-\psi)(x-h)-(u-\psi)(x)]
	\end{equation}
	then clearly $v_2\in W^{1,p(x)}_0(\Omega)$ and (\ref{eq5.2}) still trivially satisfied , due to the fact that
	\begin{equation*}
		u(x)-\psi(x)+tv_2(x)
	\end{equation*}
	\begin{equation*}
		=u(x)-\psi(x)+t\eta^2(x-h)[(u-\psi)(x-h)-(u-\psi)(x)]
	\end{equation*}
	\begin{equation*}
		=t\eta^2(x-h)(u-\psi)(x-h)+(1-t\eta^2(x-h))(u-\psi)(x)\geq 0
	\end{equation*}
	Choosing $\phi=u+tv_2$ in (\ref{eq1.3}) as test function , where $v_2$ is defined in (\ref{eq5.5}), we get
	\begin{equation}\label{eq5.6}
		\int_\Omega\big<\mathcal{A}(x,u(x),Du(x)),D(\eta^2(x-h)\tau_{-h}(u-\psi))\big>dx\geq \int_\Omega\mathcal{B}(x,u(x),Du(x)).[\eta^2(x-h)\tau_{-h}(u-\psi)]dx
	\end{equation}
	Now adding inequalities (\ref{eq5.4}) and (\ref{eq5.6}), we have
	\begin{equation*}
		\int_\Omega\big<\mathcal{A}(x,u(x),Du(x)),D(\eta^2(x)\tau_h(u-\psi))\big>dx + \int_\Omega\big<\mathcal{A}(x,u(x),Du(x)),D(\eta^2(x-h)\tau_{-h}(u-\psi))\big>dx\geq  
	\end{equation*}
	\begin{equation}\label{eq5.7}
		\int_\Omega\mathcal{B}(x,u(x),Du(x)).[\eta^2(x)\tau_h(u-\psi)+ \eta^2(x-h)\tau_{-h}(u-\psi)]dx
	\end{equation}
	By means of simple change of variable , we can write the second integral on the left hand side of the previous inequality as follows
	\begin{equation}\label{eq5.8}
		-\int_\Omega\big<\mathcal{A}(x+h,u(x+h),Du(x+h)),D(\eta^2(x)\tau_h(u-\psi))\big>dx
	\end{equation}
	and now we can write the inequality (\ref{eq5.7}) as follows
	\begin{eqnarray*}
		I=\int_\Omega\big<\mathcal{A}(x+h,u(x+h),Du(x+h))- \mathcal{A}(x,u(x),Du(x)),D(\eta^2(x)\tau_h(u-\psi))\big>dx\\
		\leq -	\int_\Omega\mathcal{B}(x,u(x),Du(x)).[\eta^2(x)\tau_h(u-\psi)+ \eta^2(x-h)\tau_{-h}(u-\psi)]dx \\
		=\int_\Omega\mathcal{B}(x,u(x),Du(x))\tau_{-h}(\eta^2\tau_h(u-\psi))dx=J
	\end{eqnarray*}
	Now we can write the integral $I$ as follows
	\begin{equation*}
		I=\int_\Omega\big<\mathcal{A}(x+h,u(x+h),Du(x+h))- \mathcal{A}(x+h,u(x+h),Du(x)),\eta^2(x)D\tau_hu\big>dx
	\end{equation*} 
	\begin{equation*}
		\hspace{8mm}-\int_\Omega\big<\mathcal{A}(x+h,u(x+h),Du(x+h))- \mathcal{A}(x+h,u(x+h),Du(x)),\eta^2(x)D\tau_h\psi\big>dx
	\end{equation*}
	\begin{equation*}
		\hspace{15mm}+\int_\Omega\big<\mathcal{A}(x+h,u(x+h),Du(x+h))- \mathcal{A}(x+h,u(x+h),Du(x)),2\eta D\eta\tau_h(u-\psi)\big>dx
	\end{equation*}
	\begin{equation*}
		+\int_\Omega\big<\mathcal{A}(x+h,u(x+h),Du(x))- \mathcal{A}(x,u(x),Du(x)),\eta^2(x)D\tau_hu\big>dx
	\end{equation*}
	\begin{equation*}
		\hspace{8mm}-\int_\Omega\big<\mathcal{A}(x+h,u(x+h),Du(x))- \mathcal{A}(x,u(x),Du(x)),\eta^2(x)D\tau_h\psi\big>dx
	\end{equation*}
	\begin{equation*}
		+\int_\Omega\big<\mathcal{A}(x+h,u(x+h),Du(x))- \mathcal{A}(x,u(x),Du(x)),2\eta D\eta\tau_h(u-\psi)\big>dx
	\end{equation*}
	\begin{equation}\label{eq5.9}
		=I_1+I_2+I_3+I_4+I_5+I_6
	\end{equation}
	that yields 
	\begin{equation}\label{eq5.10}
		I_1\leq |I_2|+|I_3|+|I_4|+|I_5|+|I_6|+|J|
	\end{equation}
	we proceed estimating the right hand side of  (\ref{eq5.10}). Please note that we denote $\mathcal{H}_i$ for $i=1,2,3,...,15$ as finite numbers depending only on  $n,\nu,L,l,\gamma_1,\gamma_2,R$ where $\mathcal{H}_j$ additionally depend on $D,\alpha$ as well for $j=4,5,6,...,12$ and for any ball $B_R\subset\Omega$ . \\[2mm]
	{\bf Estimate for $I_1$ :} the ellipticity assumption 
	(\ref{eqa1}) implies
	\begin{equation}\label{eq5.11}
		I_1\geq \nu\int_\Omega\eta^2|\tau_h Du|^2(1+|Du(x+h)|^2+|Du(x)|^2)^\frac{p(x)-2}{2}dx
	\end{equation}
	{\bf Estimate for $I_2$ :} From the growth condition (\ref{eqa2}) and Young's inequaity, we get
	\begin{equation*}
		|I_2|\leq \int_\Omega|\mathcal{A}(x+h,u(x+h),Du(x+h))-\mathcal{A}(x+h,u(x+h),Du(x))||\eta^2D\tau_h\psi|
	\end{equation*}
	\begin{equation*}
		\leq L\int_\Omega\eta^2|\tau_hDu|(1+|Du(x+h)|^2+|Du(x)|^2)^\frac{p(x)-2}{2}|D\tau_h\psi|dx
	\end{equation*}
	\begin{equation*}\footnote{Young's Inequality with $p=\frac{1}{2}$ and $q=\frac{1}{2}$}
		\leq \epsilon\int_\Omega\eta^2|\tau_hDu|^2(1+|Du(x+h)|^2+|Du(x)|^2)^\frac{p(x)-2}{2}dx+ c_\epsilon(L)\int_\Omega\eta^2|\tau_hD\psi|^2(1+|Du(x+h)|^2+|Du(x)|^2)^\frac{p(x)-2}{2}dx
	\end{equation*}
	\begin{eqnarray*}
		\leq \epsilon\int_\Omega\eta^2|\tau_hDu|^2(1+|Du(x+h)|^2+|Du(x)|^2)^\frac{p(x)-2}{2}dx + c_\epsilon(L)\Bigg(\int_{B_r}|\tau_hD\psi|^{p^+}\Bigg)^\frac{2}{p^+}\\\times \Bigg(\int_{B_r}(1+|Du(x+h)|^2+|Du(x)|^2)^{\frac{p(x)-2}{2}.\frac{p^+}{p^+-2}}dx\Bigg)^\frac{p^+-2}{p^+}
	\end{eqnarray*}
	where we used H$\ddot{o}$lder's inequality with exponent $\frac{p^+}{2}$ and $\frac{p^+}{p^+-2}$ and the properties of $\eta$.\\Here and in the sequel we denote with $\epsilon>0$ a constant that will be determined later.\\
	Observing that, since $p(x)<p^+$
	\begin{equation*}
		p^->2\implies \frac{p^+}{p^+-2}\leq \frac{p(x)}{p(x)-2}
	\end{equation*}
	we get
	\begin{equation*}
		\int_{B_r}(1+|Du(x+h)|^2+|Du(x)|^2)^{\frac{p(x)-2}{2}.\frac{p^+}{p^+-2}}dx
	\end{equation*}
	\begin{equation*}
		\leq \int_{B_r}(1+|Du(x+h)|^2+|Du(x)|^2)^\frac{p(x)}{2}dx
	\end{equation*}
	\begin{equation*}
		\leq \int_{B_r}(1+|Du(x+h)|^2+|Du(x)|^2)^\frac{p^+}{2}dx
	\end{equation*}
	\begin{equation*}
		\leq c\int_{B_r}(1+|Du(x+h)|^{p^+}+|Du(x)|^{p^+})dx
	\end{equation*}
	\begin{equation}\label{eq5.12}
		\leq c\int_{B_{2r}}(1+|Du(x)|^{p^+})dx
	\end{equation}
	Where we used the second estimate of Lemma \ref{L3.3}.\\
	Now quantity in (\ref{eq5.12}) is finite because of Lemma \ref{L5.1} for $r_0=2r$. Using this (\ref{eq5.12}) in previous estimate  of $I_2$, using the assumption $D\psi\in W^{1,\gamma_2}$ and the first estimate of Lemma \ref{L3.3} in the second integral of the right hand side of previous estimate $I_2$, we obtain
	\begin{eqnarray*}
		|I_2|\leq \epsilon\int_\Omega\eta^2|\tau_hDu|^2(1+|Du(x+h)|^2+|Du(x)|^2)^\frac{p(x)-2}{2}dx + c_\epsilon(L,n,p^{+})|h|^2\Bigg(\int_{B_r}|D^2\psi|^{p^+}\Bigg)^\frac{2}{p^+}\\\times \Bigg(\int_{B_{2r}}(1+|Du(x)|^{p^+})dx\Bigg)^\frac{p^+-2}{p^+}
	\end{eqnarray*}
	\begin{equation}\label{eq5.13}
		:= \epsilon \int_\Omega\eta^2|\tau_hDu|^2(1+|Du(x+h)|^2+|Du(x)|^2)^\frac{p(x)-2}{2}dx+|h|^2\mathcal{H}_2
	\end{equation}
	where $\mathcal{H}_2 := c_\epsilon(L,n,p^{+})\Bigg(\int_{B_r}|D^2\psi|^{p^+}\Bigg)^\frac{2}{p^+}\times \Bigg(\int_{B_{2r}}(1+|Du(x)|^{p^+})dx\Bigg)^\frac{p^+-2}{p^+}$ is finite in view of our assumptions. Arguing analogously, once more by assumption (\ref{eqa2}),I get
	\begin{equation*}
		|I_3|\leq \int_\Omega|\mathcal{A}(x+h,u(x+h),Du(x+h))-\mathcal{A}(x+h,u(x+h),Du(x))||2\eta D\eta\tau_h(u-\psi)|dx
	\end{equation*}
	\begin{equation*}
		\leq \int_\Omega 2L|\tau_hDu|(1+|Du(x+h)|^2+|Du(x)|^2)^\frac{p(x)-2}{2}\eta|D\eta||\tau_h(u-\psi)|dx
	\end{equation*}
	\begin{equation*}\footnote{Using Young's Inequality with $p,q=\frac{1}{2}$}
		\leq \epsilon\int_\Omega\eta^2|\tau_hDu|^2(1+|Du(x+h)|^2+|Du(x)|^2)^\frac{p(x)-2}{2}dx+c_\epsilon(L)\int_\Omega|D\eta|^2|\tau_h(u-\psi)|^2(1+|Du(x+h)|^2+|Du(x)|^2)^\frac{p(x)-2}{2}dx
	\end{equation*}
	\begin{equation}\label{eq5.14}
		\leq\epsilon\int_\Omega\eta^2|\tau_hDu|^2(1+|Du(x+h)|^2+|Du(x)|^2)^\frac{p(x)-2}{2}dx+\frac{c_\epsilon(L,p^+)}{r^2}\Bigg(\int_{B_r}|\tau_h(u-\psi)|^{p^+}dx\Bigg)^\frac{2}{p^+}\Bigg(\int_{B_{2r}}(1+|Du(x)|^{p^+})dx\Bigg)^\frac{p^+-2}{p^+}
	\end{equation}
	where (\ref{eq5.14}) we have used the same argument for $|I_2|$ replacing $D\psi$ with $u-\psi$ and used $|D\eta|\leq \frac{c}{r}$.\\[0.8mm]
	Since $u-\psi\in W^{1,p^+}(B_r)$, we use the first result of Lemma \ref{L3.3} to the last integral in the right hand side of the previous estimate (\ref{eq5.14}) , we will have
	\begin{eqnarray*}
		|I_3|\leq \epsilon\int_\Omega\eta^2|\tau_hDu|^2(1+|Du(x+h)|^2+|Du(x)|^2)^\frac{p(x)-2}{2}dx + |h|^2\frac{c_\epsilon(L,n,p^+)}{r^2}\Bigg(\int_{B_{2r}}|D(u-\psi)(x)|^{p^+}dx\Bigg)^\frac{2}{p^+}\\\Bigg(\int_{B_{2r}}(1+|Du(x)|^{p^+})dx\Bigg)^\frac{p^+-2}{p^+}
	\end{eqnarray*}
	\begin{equation}\label{eq5.15}
		:= \epsilon\int_\Omega\eta^2|\tau_hDu|^2(1+|Du(x+h)|^2+|Du(x)|^2)^\frac{p(x)-2}{2}+\frac{|h|^2}{r^2}\mathcal{H}_3
	\end{equation}
	where $\mathcal{H}_3 := c_\epsilon(L,n,p^+)\Bigg(\int_{B_{2r}}|D(u-\psi)(x)|^{p^+}dx\Bigg)^\frac{2}{p^+}\Bigg(\int_{B_{2r}}(1+|Du(x)|^{p^+})dx\Bigg)^\frac{p^+-2}{p^+}$ is finite.\\[0.8mm]
	Now using (\ref{eqa5}),(\ref{eq3.18}), estimate $I_4$ (see \cite{FG48_1}),
	\begin{equation*}
		|I_4|\leq \int_\Omega|\mathcal{A}(x+h,u(x+h),Du(x))-\mathcal{A}(x,u(x),Du(x))||\eta^2||\tau_hDu|dx
	\end{equation*}
	\begin{equation*}
		\leq \int_\Omega\eta^2|\tau_hDu|\;\Bigg((\kappa(x)+\kappa(x+h))+\frac{D^{-\alpha}}{e+|Du|}\Bigg)|h|(1+|Du|^2)^\frac{p(x)-1}{2}\log(e+|Du|^2)dx
	\end{equation*}
	\begin{equation*}
		\leq \int_\Omega\eta^2|h||\tau_hDu|\;\Bigg((\kappa(x)+\kappa(x+h))+\frac{D^{-\alpha}}{e+|Du|}\Bigg)(1+|Du|^2)^\frac{p(x)-1}{2}\log(e+|Du|^2)dx
	\end{equation*}
	Now we use an inequality
	\begin{equation*}
		\log(e+|Du(x)|^2)\leq \log[(e+|Du(x)|)^2]\leq 2\log (e+|Du(x)|)
	\end{equation*} 
	Now we have by Young's Inequality
	\begin{eqnarray*}\footnote{Using Young's Inequality with $p,q=\frac{1}{2}$}
		\leq \epsilon\int_\Omega\eta^2|\tau_hDu|^2\;(1+|Du(x)|^2+|Du(x+h)|^2)^\frac{p(x)-2}{2}dx+c_\epsilon \int_{B_r}|h|^2 \Bigg((\kappa(x)+\kappa(x+h))+\frac{D^{-\alpha}}{e+|Du|}\Bigg)^2\\(1+|Du(x)|^2)^\frac{p(x)}{2}\log^2(e+|Du(x)|)dx
	\end{eqnarray*}
	\begin{eqnarray*}
		\leq \epsilon\int_\Omega\eta^2|\tau_hDu|^2\;(1+|Du(x)|^2+|Du(x+h)|^2)^\frac{p(x)-2}{2}dx+c_\epsilon |h|^2\int_{B_r} \Bigg((\kappa(x)+\kappa(x+h))+\frac{D^{-\alpha}}{e+|Du|}\Bigg)^2\\(1+|Du(x)|^2)^\frac{p(x)}{2}\log^2(e+|Du(x)|)dx
	\end{eqnarray*}
	\begin{eqnarray*}
		\leq \epsilon\int_\Omega\eta^2|\tau_hDu|^2\;(1+|Du(x)|^2+|Du(x+h)|^2)^\frac{p(x)-2}{2}dx+c_\epsilon |h|^2\int_{B_r} \Bigg((\kappa(x)+\kappa(x+h))+\frac{D^{-\alpha}}{e+|Du|}\Bigg)^2\\(1+|Du(x)|^2)^\frac{p(x)}{2}\log^2(e+|Du(x)|)dx
	\end{eqnarray*}
	\begin{eqnarray*}
		\leq \epsilon\int_\Omega\eta^2|\tau_hDu|^2\;(1+|Du(x)|^2+|Du(x+h)|^2)^\frac{p(x)-2}{2}dx+c_\epsilon|h|^2 \Bigg( \int_{B_r} \Bigg((\kappa(x)+\kappa(x+h))+\frac{D^{-\alpha}}{e+|Du|}\Bigg)^n\\\log^n(e+|Du(x)|)dx \Bigg)^\frac{2}{n}\times\Bigg(\int_{B_r}(1+|Du(x)|^2)^\frac{np(x)}{2(n-2)}dx\Bigg)^\frac{n-2}{n}
	\end{eqnarray*}
\begin{eqnarray*}
	\leq \epsilon\int_\Omega\eta^2|\tau_hDu|^2\;(1+|Du(x)|^2+|Du(x+h)|^2)^\frac{p(x)-2}{2}dx+c_\epsilon|h|^2 \Bigg( \int_{B_r} (\kappa^n(x+h)+\kappa^n(x))\log^n(e+|Du(x)|)dx \\+\int_{B_r}\frac{D^{-n\alpha}}{(e+|Du|)^n}\log^n(e+|Du(x)|)dx\Bigg)^\frac{2}{n}\times\Bigg(\int_{B_r}(1+|Du(x)|^2)^\frac{np(x)}{2(n-2)}dx\Bigg)^\frac{n-2}{n}
\end{eqnarray*}
\begin{eqnarray*}
	\leq \epsilon\int_\Omega\eta^2|\tau_hDu|^2\;(1+|Du(x)|^2+|Du(x+h)|^2)^\frac{p(x)-2}{2}+c_\epsilon|h|^2 \Bigg( \int_{B_r} (\kappa^n(x+h)+\kappa^n(x))\log^n(e+|Du(x)|)dx \\+\int_{B_r}\frac{D^{-n\alpha}}{(e+|Du|)^n}(e+|Du(x)|)^ndx\Bigg)^\frac{2}{n}\times\Bigg(\int_{B_r}(1+|Du(x)|^2)^\frac{np(x)}{2(n-2)}dx\Bigg)^\frac{n-2}{n}
\end{eqnarray*}
	\begin{eqnarray*}
		\leq \epsilon\int_\Omega\eta^2|\tau_hDu|^2\;(1+|Du(x)|^2+|Du(x+h)|^2)^\frac{p(x)-2}{2}+c_\epsilon|h|^2 \Bigg( \int_{B_r} (\kappa^n(x+h)+\kappa^n(x))\log^n(e+|Du(x)|)dx +r^n\Bigg)^\frac{2}{n}\\\times\Bigg(\int_{B_r}(1+|Du(x)|^2)^\frac{np(x)}{2(n-2)}dx\Bigg)^\frac{n-2}{n}
	\end{eqnarray*}
	
	Now consider the term 
	\begin{equation*}
		\int_{B_r} (\kappa^n(x+h)+\kappa^n(x))\log^n(e+|Du(x)|)dx
	\end{equation*} 
	\begin{equation}\label{eq5.16}
		\leq \int_{B_r} \kappa^n(x+h)\log^n(e+|Du(x)|)dx+\int_{B_r} \kappa^n(x)\log^n(e+|Du(x)|)dx
	\end{equation}
	Now put $\epsilon=1,\gamma=1,\alpha=n,s=\kappa^n(x+h),t=\log^n(e+|Du(x)|)$, then by Lemma \ref{L3.2}
	\begin{equation*}
		\int_{B_r} \kappa^n(x+h)\log^n(e+|Du(x)|)dx
	\end{equation*}
	\begin{equation*}
		\leq \int_{B_r}\kappa^n(x+h)\log^n(e+\kappa^n(x+h))dx+\int_{B_r}(e+|Du(x)|-1)\log^n(e+|Du(x)|)dx
	\end{equation*}
	\begin{equation}\label{eq5.17}
		\leq c \int_{B_{2r}}\kappa^n(x+h)\log^n(e+\kappa(x+h))+\int_{B_r}(e+|Du(x)|-1)\log^n(e+|Du(x)|)dx
	\end{equation}
	Again similarly by Lemma \ref{L3.2}, we have for the second integral of the right hand side of (\ref{eq5.16}) with $\epsilon=1,\gamma=1,\alpha=n,s=\kappa^n(x),t=\log^n(e+|Du(x)|)$ we get
	\begin{equation*}
		\int_{B_r} \kappa^n(x)\log^n(e+|Du(x)|)dx
	\end{equation*}
	\begin{equation*}
		\leq \int_{B_r}\kappa^n(x)\log^n(e+\kappa^n(x))dx+\int_{B_r}(e+|Du(x)|-1)\log^n(e+|Du(x)|)dx
	\end{equation*}
	\begin{equation}\label{eq5.18}
		\leq c \int_{B_{2r}}\kappa^n(x)\log^n(e+\kappa(x))+\int_{B_r}(e+|Du(x)|-1)\log^n(e+|Du(x)|)dx
	\end{equation}
	Now from (\ref{eq5.17}) and (\ref{eq5.18}), we have from(\ref{eq5.16})
	\begin{equation*}
		\int_{B_r} (\kappa^n(x+h)+\kappa^n(x))\log^n(e+|Du(x)|)dx
	\end{equation*} 
	\begin{equation}\label{eq5.19}
		\leq c \int_{B_{2r}}\kappa^n(x)\log^n(e+\kappa(x))+\int_{B_r}(e+|Du(x)|-1)\log^n(e+|Du(x)|)dx
	\end{equation}
	Due to assumption of $\kappa\in L^n\log^nL(\Omega)$, the first integral to the right hand side of (\ref{eq5.19}) is finite. Now we want to recall an inequality for any $a,b>0$,
	\begin{equation}\label{eq5.20}
		\log(e+ab)\leq \log(e+a)+\log(e+b)
	\end{equation}
	Now consider second integral of the right hand side of (\ref{eq5.19})
	\begin{equation*}
		\int_{B_r}(e-1+|Du(x)|)\log^n(e+|Du(x)|)dx
	\end{equation*}
	\begin{equation*}
		\leq c\int_{B_r\cap\{|Du|\geq e\}}|Du(x)||\log^n(e+|Du(x)|)|dx+\int_{B_r\cap\{|Du|<e\}}(e+|Du(x)|)\log^n(e+|Du(x)|)
	\end{equation*}
	\begin{equation*}
		\leq c|B_r|-\!\!\!\!\!\!\!\int_{B_r}|Du(x)||\log^n\Bigg(e+\frac{|Du(x)}{||\;|Du|\;||_1}\Bigg)dx + c|B_r|-\!\!\!\!\!\!\!\int_{B_r}|Du(x)|\log^n(e+||\;|Du|\;||_1)dx+c|B_r|
	\end{equation*}
	\begin{equation}\label{eq5.21}
		\leq c|B_r|\Bigg(-\!\!\!\!\!\!\!\int_{B_r}|Du(x)|^2dx\Bigg)^\frac{1}{2}+c|B_r|-\!\!\!\!\!\!\!\int_{B_r}|Du(x)|\log^n(e+||\;|Du|\;||_1)dx+c|B_r|
	\end{equation}
	First term of (\ref{eq5.21}) is bounded because of $u\in W^{1,p(x)}(\Omega)$ with $p(x)>2$. Now we have from second integral (see for reference \cite[p. 132]{AM_56_96_6_2})
	\begin{equation*}
		c|B_r|-\!\!\!\!\!\!\!\int_{B_r}|Du(x)|\log^n(e+||\;|Du|\;||_1)dx\leq c\log^n\Bigg(\frac{1}{r}\Bigg)\int_{B_r}|Du(x)|dx\leq \frac{c}{r^n}\int_{B_r}|Du(x)|dx
	\end{equation*}
	Since in Lemma \ref{L5.1} , we can have $Du\in L^{p^+}(\Omega)$ and $p^+>1$ so by Lebesgue space inclusion for bounded domain , we have $Du\in L^1(\Omega)$ so second term is finite precisely $\int_{B_r}|Du(x)|dx<
	+\infty$ \vskip 2mm
	Now aim to estimate 
	\begin{equation*}
		\Bigg(\int_{B_r}(1+|Du(x)|^2)^\frac{np(x)}{2(n-2)}\Bigg)^\frac{n-2}{n}
	\end{equation*}
	This is the crucial point where the Calder\'on-Zygmund result is used.\\[0.8mm]
	Since $D\psi\in W^{1,p^+}(\Omega)$, classical Sobolev embedding theorem implies $D\psi\in L^{(p^+)^*}$. Observing that
	\begin{equation*}\footnote{$\frac{1}{p^*}=\frac{1}{p}-\frac{1}{n}$ for $p<n$}
		p^+\geq 2\implies\frac{np(x)}{n-2}\leq \frac{np^+}{n-p^+}=(p^+)^*
	\end{equation*}
	I have $|D\psi|^{p(x)}\in L^{\frac{n}{n-2}}_{loc}(\Omega)$. Since $\frac{\gamma_1}{\gamma_1-2}>\frac{n}{n-2}$, then by assumption on $q$ , $a^\frac{p(x)}{p(x)-1}\in  L^{\frac{n}{n-2}}_{loc}(\Omega)$ \\[0.8mm]
	Therefore, by Calder\'on- Zygmund Theorem with $|Du|^{p(x)}\in L^\frac{n}{n-2}(\Omega)$, in perticular from (\ref{eq2.8}) there holds
	\begin{equation}\label{eq5.22}
		\Bigg[\int_{B_r}|Du|^\frac{np(x)}{n-2}dx\Bigg]^\frac{n-2}{n}\leq c\bigg[\int_{B_r}a^{\frac{p(x)}{p(x)-1}.\frac{n}{n-2}}+|u|^\frac{np(x)}{n-2}+|D\psi|^\frac{np(x)}{n-2}+|Du|^{p(x)}+1\bigg]^\frac{n-2}{n}dx
	\end{equation}
	where $c=c(n,\nu',\gamma_1,\gamma_2,L',q)$. Since all the term on the right hand side of (\ref{eq5.22}) are finite so we have
	\begin{equation}\label{eq5.23}
		\bigg(\int_{B_r}\Big(1+|Du(x)|\Big)^\frac{np(x)}{(n-2)}\bigg)^\frac{n-2}{n}\leq cH'
	\end{equation}
	So I have 
	\begin{eqnarray*}
		|I_4|\leq \epsilon\int_\Omega\eta^2|\tau_hDu|^2\;(1+|Du(x)|^2+|Du(x+h)|^2)^\frac{p(x)-2}{2}dx+ c|h|^2\Bigg[\int_{\Omega}\kappa^n(x)\log^n(e+\kappa(x))dx+\Bigg(-\!\!\!\!\!\!\!\int_{B_r}|Du(x)|^2dx\Bigg)^\frac{1}{2}\\+\frac{1}{r^n}\int_{B_r}|Du(x)|dx+2r^n\Bigg]^\frac{2}{n}\\\times CH'
	\end{eqnarray*}
	\begin{equation*}
		\leq \epsilon\int_\Omega\eta^2|\tau_hDu|^2\;(1+|Du(x)|^2+|Du(x+h)|^2)^\frac{p(x)-2}{2}dx+ \frac{c|h|^2}{r^2}\Bigg[\mathcal{H}_4+\mathcal{H}_5r^{2n}+\mathcal{H}_6r^n\Bigg]^\frac{2}{n}
	\end{equation*}
	\begin{equation*}
		\leq \epsilon\int_\Omega\eta^2|\tau_hDu|^2\;(1+|Du(x)|^2+|Du(x+h)|^2)^\frac{p(x)-2}{2}dx+ \frac{c|h|^2}{r^2}\Bigg[\mathcal{H}_4+\mathcal{H}_5r^{2n}+\mathcal{H}_6r^n+1\Bigg]
	\end{equation*}
	\begin{equation}\label{eq5.24}
		:= \epsilon\int_\Omega\eta^2|\tau_hDu|^2\;(1+|Du(x)|^2+|Du(x+h)|^2)^\frac{p(x)-2}{2}dx+ \Bigg[\mathcal{H}_4\frac{|h|^2}{r^2}+\mathcal{H}_5|h|^2r^{2n-2}+\mathcal{H}_6|h|^2r^{n-2}\Bigg]
	\end{equation}
	{\bf Estimate for $I_5$: } Arguing as we did for the estimate of the integral $I_4$, we get with help of (\ref{eqa5}),(\ref{eq3.18})  and H$\ddot{o}$lder's Inequlaity
	\begin{equation*}
		|I_5|\leq \int_\Omega|\eta|^2|\mathcal{A}(x+h,u(x+h),Du(x))-\mathcal{A}(x,u(x),Du(x))||\tau_hD\psi|
	\end{equation*}
	\begin{equation*}
		\leq \int_\Omega|\eta|^2(1+|\tau_hu|)|h|\Bigg((\kappa(x)+\kappa(x+h))+\frac{D^{-\alpha}}{e+|Du|}\Bigg)(1+|Du(x)|^2)^\frac{p(x)-1}{2}\log(e+|Du(x)|^2)|\tau_hD\psi|dx
	\end{equation*}
	\begin{equation*}
		\leq |h|\int_\Omega|\eta|^2\Bigg((\kappa(x)+\kappa(x+h))+\frac{D^{-\alpha}}{e+|Du|}\Bigg)(1+|Du(x)|^2)^\frac{p(x)-1}{2}\log(e+|Du(x)|^2)|\tau_hD\psi|dx
	\end{equation*}
	\begin{equation*}
		\leq c|h|\Bigg(\int_{B_r}|\tau_hD\psi|^{p^+}dx\Bigg)^\frac{1}{p^+}\Bigg(\int_{B_r}\Bigg(\Big((\kappa(x)+\kappa(x+h))+\frac{D^{-\alpha}}{e+|Du|}\Big)\log(e+|Du(x)|^2)dx\Bigg)^\frac{p^+}{p^+-1}\Big(1+|Du(x)|^2\Big)^{\frac{p(x)-1}{2}.\frac{p^+}{p^+-1}}dx\Bigg)^\frac{p^+-1}{p^+}
	\end{equation*}
	\begin{equation*}\footnote{$\frac{p^+}{p^+-1}\leq \frac{p(x)}{p(x)-1}$}
		\leq c|h|\Bigg(\int_{B_r}|\tau_hD\psi|^{p^+}dx\Bigg)^\frac{1}{p^+}\Bigg(\int_{B_r}\Bigg(\Big((\kappa(x)+\kappa(x+h))+\frac{D^{-\alpha}}{e+|Du|}\Big)\log(e+|Du(x)|^2)dx\Bigg)^\frac{p^+}{p^+-1}\Big(1+|Du(x)|^2\Big)^\frac{p(x)}{2}dx\Bigg)^\frac{p^+-1}{p^+}
	\end{equation*}
	\begin{equation*}
		\leq c|h|\Bigg(\int_{B_r}|\tau_hD\psi|^{p^+}dx\Bigg)^\frac{1}{p^+}\Bigg(\int_{B_r}\Bigg(\Big((\kappa(x)+\kappa(x+h)\Big)+\frac{D^{-\alpha}}{e+|Du|}\Bigg)^n\log^n(e+|Du(x)|^2)dx\Bigg)^\frac{1}{n}
	\end{equation*}
\begin{equation}\label{eq5.25}
	\Bigg(\int_{B_r}\Big(1+|Du(x)|^2\Big)^{\frac{p(x)}{2}.\frac{n(p^+-1)}{n(p^+-1)-2}}dx\Bigg)^\frac{p^+-1}{p^+}
\end{equation}
	Since , $|D\psi|\in L^{(p^+)^*}(\Omega)$. Since
	\begin{equation*}
		\frac{np(x)(p^+-1)}{n(p^+-1)-p^+}\leq \frac{np^+}{n-p^+}=(p^+)^*\iff p^+\geq 2
	\end{equation*}
	we have $|D\psi|^{p(x)}\in L^\frac{n(p^+-1)}{n(p^+-1)-p^+}(\Omega)$. Since $\frac{n}{n-\frac{\gamma_1}{\gamma_1-2}}>\frac{n(p^+-1)}{n(p^+-1)-p^+}$ , so by assumption on $q$, $a(x)^\frac{p(x)}{p(x)-1}\in L^\frac{n(p^+-1)}{n(p^+-1)-p^+}(\Omega)$ \\
	So by Calder\'on- Zygmund Estimate (\ref{eq2.8})$, |Du|^{p(x)}\in L^{\frac{n(p^+-1)}{n(p^+-1)-p^+}}(\Omega)$\\ 
	Do by the help of Calder\'on-Zygmund eatimates (\ref{eq2.8}) we have
	\begin{eqnarray*}
		\bigg(\int_{B_r}\Big(1+|Du(x)|\Big)^\frac{np(x).(p^+-1)}{n(p^+-1)-p^+}\bigg)^\frac{n(p^+-1)-p^+}{n(p^+-1)}\leq c\bigg[\int_{B_r}\Big(a^{\frac{p(x)}{p(x)-1}.\frac{n(p^+-1)}{n(p^+-1)-p^+}}+|u|^{p(x).\frac{n(p^+-1)}{n(p^+-1)-p^+}}+|D\psi|^{p(x).\frac{n(p^+-1)}{n(p^+-1)-p^+}}+|Du|^{p(x)}\\1\Big)dx\bigg]^\frac{n(p^+-1)-p^+}{n(p^+-1)}
	\end{eqnarray*}
	\begin{equation}\label{eq5.26}	\bigg(\int_{B_r}\Big(1+|Du(x)|\Big)^\frac{np(x).(p^+-1)}{n(p^+-1)-p^+}\bigg)^\frac{n(p^+-1)-p^+}{n(p^+-1)}\leq cH''
	\end{equation}
	Inserting (\ref{eq5.26}) into (\ref{eq5.25}), we have by the first part of Lemma \ref{L3.3} and $I_4$ estimate
	\begin{equation*}
		|I_5|\leq c|h|\Bigg(\int_{B_r}|\tau_hD\psi|^{p^+}dx\Bigg)^\frac{1}{p^+}\Bigg(\int_{B_r}\Bigg(\Big((\kappa(x)+\kappa(x+h)\Big)+\frac{D^{-\alpha}}{e+|Du|}\Bigg)^n\log^n(e+|Du(x)|^2)dx\Bigg)^\frac{1}{n}CH''
	\end{equation*}
	\begin{equation*}\footnote{The underbrace term is finite due to $D\psi\in W^{1,p^+}(\Omega)$}
		\leq c(n,p^+)|h|^2\underbrace{\Bigg(\int_{B_r}|D^2\psi|^{p^+}dx\Bigg)^\frac{1}{p^+}}\Bigg(\int_{B_r}\Bigg(\Big((\kappa(x)+\kappa(x+h)\Big)+\frac{D^{-\alpha}}{e+|Du|}\Bigg)^n\log^n(e+|Du(x)|^2)dx\Bigg)^\frac{1}{n} CH''
	\end{equation*}
Then by estiamte $I_4$, we have
	\begin{equation*}
		\leq |h|^2\Big(\mathcal{H}_7+\mathcal{H}_8r^n+\frac{\mathcal{H}_9}{r^n}\Big)^\frac{1}{n}
	\end{equation*}
	\begin{equation}\label{eq5.27}
		:= \mathcal{H}_7|h|^2r^{n-1}+\mathcal{H}_8|h|^2r^{2n-1}+\mathcal{H}_9\frac{|h|^2}{r}
	\end{equation}
	{\bf Estiamte of $I_6$ :} Using (\ref{eqa5}), H$\ddot{o}$lder's Inequality and (\ref{eq3.18})
	\begin{equation*}
		|I_6|\leq \int_\Omega|\mathcal{A}(x+h,u(x+h),Du(x))-\mathcal{A}(x,u(x),Du(x))||2\eta D\eta\tau_h(u-\psi)|dx
	\end{equation*}
	\begin{equation*}
		\leq 2|h|\int_\Omega\eta|D\eta||\tau_h(u-\psi)(x)|\Bigg((\kappa(x)+\kappa(x+h))+\frac{D^{-\alpha}}{e+|Du|}\Bigg)(1+|Du(x)|^2)^\frac{p(x)-1}{2}\log(e+|Du(x)|^2)dx
	\end{equation*}
	\begin{eqnarray*}\footnote{Using Holder's Inequality and $|D\eta|\leq \frac{c}{r}$ in $B_r$}
		\leq \frac{c}{r}|h|\Bigg(\int_{B_r}|\tau_h(u-\psi)(x)|^{p^+}dx\Bigg)^\frac{1}{p^+}\Bigg(\int_{B_r}\Bigg(\Big((\kappa(x)+\kappa(x+h))+\frac{D^{-\alpha}}{e+|Du|}\Big)\Bigg)^\frac{p^+}{p^+-1}(\log(e+|Du(x)|))^\frac{p^+}{p^+-1}\\(1+|Du(x)|^2)^\frac{p(x)}{2}dx\Bigg)^\frac{p^+-1}{p^+}
	\end{eqnarray*}
	\begin{eqnarray*}\footnote{Using Holder's inequality}
		\leq \frac{c}{r}|h|\Bigg(\int_{B_r}|\tau_h(u-\psi)(x)|^{p^+}\Bigg)^\frac{1}{p^+}\Bigg(\int_{B_r}\Bigg((\kappa(x)+\kappa(x+h))+\frac{D^{-\alpha}}{e+|Du|}\Bigg)^n\log^n(e+|Du(x)|)dx\Bigg)^\frac{1}{n}\\\Bigg(\int_{B_r}(1+|Du(x)|)^\frac{p(x)n(p^+-1)}{n(p^+-1)-p^+}dx\Bigg)^\frac{n(p^+-1)-p^+}{np^+}
	\end{eqnarray*}
	\begin{equation*}
		\leq c(n,p^+)\frac{|h|^2}{r}\Bigg(\int_{B_{2r}}|D(u-\psi)(x)|^{p^+}\Bigg)^\frac{1}{p^+}\Bigg(\int_{B_r}\Bigg((\kappa(x)+\kappa(x+h))^n+\frac{D^{-n\alpha}}{(e+|Du|)^n}\Bigg)\log^n(e+|Du(x)|)dx\Bigg)^\frac{1}{n}
	\end{equation*}
\begin{equation}\label{eq5.28}
	\Bigg(\int_{B_r}(1+|Du(x)|)^\frac{p(x)n(p^+-1)}{n(p^+-1)-p^+}dx\Bigg)^\frac{n(p^+-1)-p^+}{np^+}
\end{equation}
	In the first integral of last line we use first result of Lemma \ref{L3.3}. Now since by assumption $\psi\in W^{1,\gamma_2}(\Omega)$ and $\gamma_2\geq p^+$, so $D\psi\in L^{p^+}(\Omega)$. From Lemma \ref{L5.1}, $Du\in L^{p^+}(\Omega)$, so clearly $D(u-\psi)\in L^{p^+}(\Omega)$. Now with the help of $I_4$ estimate and (\ref{eq5.26}), we can conclude
	\begin{equation*}
		\leq \frac{|h|^2}{r}\Big(\mathcal{H}_{10}+\mathcal{H}_{11}r^n+\mathcal{H}_{12}\frac{1}{r^n}\Big)^\frac{1}{n}
	\end{equation*}
	\begin{equation}\label{eq5.29}
		:= \mathcal{H}_{10}|h|^2r^{n-2}+\mathcal{H}_{11}|h|^2r^{2n-2}+\mathcal{H}_{12}\frac{|h|^2}{r^2}
	\end{equation}
	{\bf Estimate for $J$:} With help of assumption (\ref{eqa4}),\\
	\begin{equation*}
		|J|\leq \int_\Omega|\mathcal{B}(x,u,Du)||\tau_{-h}(\eta^2\tau_h(u-\psi))|dx
	\end{equation*}
	\begin{equation*}\footnote{Holder's Inequality}
		\leq \Bigg(\int_{B_r}|\mathcal{B}|^\frac{p(x)}{p(x)-1}dx\Bigg)^\frac{p(x)-1}{p(x)}\Bigg(\int_\Omega|\tau_{-h}\Big(\eta^2\tau_h(u-\psi)\Big)|^{p(x)}dx\Bigg)^\frac{1}{p(x)}
	\end{equation*}
	\begin{equation*}
		\leq c \Bigg(\int_{B_r}\Big(|Du|^\frac{r(x)p(x)}{p(x)-1}+|u|^\frac{r(x)p(x)}{p(x)-1}+|a|^\frac{p(x)}{p(x)-1}\Big)dx\Bigg)^\frac{p(x)-1}{p(x)}\Bigg(\int_\Omega|\eta^2(x-h)\tau_h(u-\psi)(x-h)-\eta^2(x)\tau_h(u-\psi)(x)|^{p(x)}dx\Bigg)^\frac{1}{p(x)}
	\end{equation*}
	\begin{equation}\label{eq5.30}
		\leq c \Bigg(\int_{B_r}\Big(|Du|^{p(x)}+|u|^{p(x)}+|a|^\frac{p(x)}{p(x)-1}+1\Big)dx\Bigg)^\frac{p(x)-1}{p(x)}\Bigg(\int_{B_r}\Big(|\eta^2(x-h)\tau_h(u-\psi)(x-h)|+|\eta^2(x)\tau_h(u-\psi)(x)|\Big)^{p(x)}dx\Bigg)^\frac{1}{p(x)}
	\end{equation}
	Now consider the second integral of (\ref{eq5.30})
	\begin{equation*}
		\Bigg(\int_{B_r}\Big(|\eta^2(x-h)\tau_h(u-\psi)(x-h)|+|\eta^2(x)\tau_h(u-\psi)(x)|\Big)^{p(x)}dx\Bigg)^\frac{1}{p(x)}
	\end{equation*}
	\begin{equation*}
		\leq c\Bigg(\int_{B_r}\Big(|\eta^2(x-h)\tau_h(u-\psi)(x-h)|+|\eta^2(x)\tau_h(u-\psi)(x)|+1\Big)^{p^+}dx+1\Bigg)
	\end{equation*}
	\begin{equation*}
		\leq c\Bigg(\int_{B_r}\Big(|\eta^2(x-h)\tau_h(u-\psi)(x-h)|^{p^+}+|\eta^2(x)\tau_h(u-\psi)(x)|^{p^+}+1\Big)dx+1\Bigg)
	\end{equation*}
	\begin{equation*}
		\leq c\Bigg(\int_{B_r}|\eta^2(x-h)\tau_h(u-\psi)(x-h)|^{p^+}dx+\int_{B_r}|\eta^2(x)\tau_h(u-\psi)(x)|^{p^+}dx+r^n+1\Bigg)
	\end{equation*}
	\begin{equation*}
		\leq c\Bigg(\int_{B_{2r}}|\eta^2(x)\tau_h(u-\psi)(x)|^{p^+}dx+\int_{B_r}|\eta^2(x)\tau_h(u-\psi)(x)|^{p^+}dx+r^n+1\Bigg)
	\end{equation*}
	\begin{equation*}
		\leq c \Bigg(\int_{B_{2r}}|\eta^2(x)\tau_h(u-\psi)(x)|^{p^+}dx+r^n+1\Bigg)
	\end{equation*}
	\begin{equation}\label{eq5.31}
		\leq c(n,p^+)\Bigg(|h|^{p^+}\int_{B_{3r}}|D(u-\psi)|^{p^+}dx+r^n+1\Bigg)
	\end{equation}
	In (\ref{eq5.31}), using first result of Lemma \ref{L3.3}. Since $u-\psi\in W^{1,p^+}_{loc}(\Omega)$ and $p^+>\gamma_1$ , then we have from (\ref{eq5.31})
	\begin{equation}\label{eq5.32}
		\leq \Bigg(\mathcal{H}_{13}|h|^{\gamma_1}+\mathcal{H}_{14}r^n+\mathcal{H}_{15}\Bigg)
	\end{equation}
	Now consider the first integral of (\ref{eq5.30})
	\begin{equation*}
		\Bigg(\int_{B_r}\Big(|Du|^{p(x)}+|u|^{p(x)}+|a|^\frac{p(x)}{p(x)-1}+1\Big)dx\Bigg)^\frac{p(x)-1}{p(x)}
	\end{equation*}
	\begin{equation*}
		\leq  \Bigg(\int_{B_r}\Big(|Du|^\frac{{p(x)n}}{n-2}+|u|^\frac{{p(x)n}}{n-2}+|a|^{\frac{p(x)n}{(p(x)-1)(n-2)}}+1\Big)dx\Bigg)^\frac{(p(x)-1)(n-2)}{n p(x)}.|r|^\frac{2(p(x)-1)}{p(x)}
	\end{equation*}
	According to (\ref{eq5.23}) , $q>\frac{n}{n-2}$, $\frac{p(x)-1}{p(x)}\geq \frac{\gamma_1-1}{\gamma_1}$ and $r<1$ , we have
	\begin{equation*}
		\leq c |r|^\frac{2(\gamma_1-1)}{\gamma_1}
	\end{equation*}
	\begin{equation}\label{eq5.33}
		\leq c |r|^\frac{2(\gamma_1-1)}{\gamma_1}
	\end{equation}
	Now inserting (\ref{eq5.32}) and (\ref{eq5.33})into (\ref{eq5.30}) implies,
	\begin{equation}\label{eq5.34}
		|J|\leq \mathcal{H}_{13}|h|^{\gamma_1}.|r|^\frac{\gamma_1-1}{\gamma_1}+\mathcal{H}_{14}|r|^{\frac{\gamma_1-1}{\gamma_1}+n}+\mathcal{H}_{15}|r|^\frac{2(\gamma_1-1)}{\gamma_1}
	\end{equation}
	Finally inserting (\ref{eq5.11},\ref{eq5.13},\ref{eq5.15},\ref{eq5.24},\ref{eq5.27},\ref{eq5.29},\ref{eq5.34}) in (\ref{eq5.10})
	\begin{eqnarray*}
		\nu\int\eta^2|\tau_hDu|^2(1+|Du(x+h)|^2+|Du(x)|^2)^\frac{p(x)-2}{2}dx\leq 3\epsilon \int\eta^2|\tau_hDu|^2(1+|Du(x+h)|^2+|Du(x)|^2)^\frac{p(x)-2}{2}dx\\+|h|^2(\mathcal{H}_2+\mathcal{H}_3)+\Bigg[\mathcal{H}_4\frac{|h|^2}{r^2}+\mathcal{H}_5|h|^2r^{2n-2}+\mathcal{H}_6|h|^2r^{n-2}\Bigg]+\Bigg[\mathcal{H}_7|h|^2r^{n-1}+\mathcal{H}_8|h|^2r^{2n-1}+\mathcal{H}_9\frac{|h|^2}{r}\Bigg]\\+\Bigg[\mathcal{H}_{10}|h|^2r^{n-2}+\mathcal{H}_{11}|h|^2r^{2n-2}+\mathcal{H}_{12}\frac{|h|^2}{r^2}\Bigg]+\Bigg[\mathcal{H}_{13}|h|^{\gamma_1}.|r|^\frac{2(\gamma_1-1)}{\gamma_1}+\mathcal{H}_{14}|r|^{\frac{2(\gamma_1-1)}{\gamma_1}+n}+\mathcal{H}_{15}|r|^\frac{2(\gamma_1-1)}{\gamma_1}\Bigg]
	\end{eqnarray*}
	Choosing $\epsilon=\frac{\nu}{6}$ I get, 
	\begin{eqnarray*}
		\nu\int\eta^2|\tau_hDu|^2(1+|Du(x+h)|^2+|Du(x)|^2)^\frac{p(x)-2}{2}dx\leq |h|^2(\mathcal{H}_2+\mathcal{H}_3)+\Bigg[\mathcal{H}_4\frac{|h|^2}{r^2}+\mathcal{H}_5|h|^2r^{2n-2}+\mathcal{H}_6|h|^2r^{n-2}\Bigg]\\+\Bigg[\mathcal{H}_7|h|^2r^{n-1}+\mathcal{H}_8|h|^2r^{2n-1}+\mathcal{H}_9\frac{|h|^2}{r}\Bigg]+\Bigg[\mathcal{H}_{10}|h|^2r^{n-2}+\mathcal{H}_{11}|h|^2r^{2n-2}+\mathcal{H}_{12}\frac{|h|^2}{r^2}\Bigg]\\+\Bigg[\mathcal{H}_{13}|h|^{\gamma_1}.|r|^\frac{2(\gamma_1-1)}{\gamma_1}+\mathcal{H}_{14}|r|^{\frac{2(\gamma_1-1)}{\gamma_1}+n}+\mathcal{H}_{15}|r|^\frac{2(\gamma_1-1)}{\gamma_1}\Bigg]
	\end{eqnarray*}
	Using Lemma \ref{L3.1} in the left hand side of previous estimate and recalling that $\eta\equiv 1$ on $B_\frac{r}{2}$, we get
	\begin{equation*}
		\nu\int\eta^2|\tau_hDu|^2(1+|Du(x+h)|^2+|Du(x)|^2)^\frac{p(x)-2}{2}dx
	\end{equation*}
	\begin{equation*}
		\geq \nu\int\eta^2|\tau_hDu|^2(1+|Du(x+h)|^2+|Du(x)|^2)^\frac{\gamma_1-2}{2}dx
	\end{equation*}
	\begin{equation*}
		\geq \nu \int_{B_{r/2}}\eta^2|\tau_hV\gamma_1(Du(x))|^2dx
	\end{equation*}
	\begin{equation*}
		\geq \nu \int_{B_{r/2}}|\tau_hV\gamma_1(Du(x))|^2dx
	\end{equation*}
	Then I have 
	\begin{eqnarray*}
		\nu \int_{B_{r/2}}|\tau_hV\gamma_1(Du(x))|^2dx\leq|h|^2(\mathcal{H}_2+\mathcal{H}_3)+\Bigg[\mathcal{H}_4\frac{|h|^2}{r^2}+\mathcal{H}_5|h|^2r^{2n-2}+\mathcal{H}_6|h|^2r^{n-2}\Bigg]\\+\Bigg[\mathcal{H}_7|h|^2r^{n-1}+\mathcal{H}_8|h|^2r^{2n-1}+\mathcal{H}_9\frac{|h|^2}{r}\Bigg]+\Bigg[\mathcal{H}_{10}|h|^2r^{n-2}+\mathcal{H}_{11}|h|^2r^{2n-2}+\mathcal{H}_{12}\frac{|h|^2}{r^2}\Bigg]
	\end{eqnarray*}
	\begin{equation}\label{eq5.35}
		\hspace{80mm}+\Bigg[\mathcal{H}_{13}|h|^{\gamma_1}.|r|^\frac{2(\gamma_1-1)}{\gamma_1}+\mathcal{H}_{14}|r|^{\frac{2(\gamma_1-1)}{\gamma_1}+n}+\mathcal{H}_{15}|r|^\frac{2(\gamma_1-1)}{\gamma_1}\Bigg]
	\end{equation}
	for every ball $B_r\subset B_{2r}\Subset \Omega$.\vskip 0.5mm
	Now fix arbitrary open sets $\Omega'\Subset\Omega''\Subset\Omega$ and choose a fixed number $y_0\in \Omega'$. Let consider $\sigma\in (0,1)$ which will be choosen later according to the situation and consider the ball $\mathfrak{B}=\mathfrak{B}(y_0,|h|^\sigma)$ with sufficiently small $|h|$, depending on the dimension $n$. we have choosen parameter $\sigma$ and the distance between $\Omega'$ and the boundary of $\Omega''$ such that $\mathcal{Q}_{out}(\hat{\mathfrak{B}})\Subset\Omega''$.\vskip 0.5mm
	Now recent estimate (\ref{eq5.35}) applied over ball $\mathfrak{B}$ implies 
	\begin{eqnarray*}
		\int_\mathfrak{B}|\tau_hV_{\gamma_1}(Du)|^2dx\leq |h|^2G_1+\Bigg[G_2|h|^{2-2\sigma}+G_3|h|^{2+2\sigma n-2\sigma}+G_4|h|^{2+\sigma n-2\sigma}\Bigg]+\Bigg[G_5|h|^{2+\sigma n-\sigma}+G_6|h|^{2+2\sigma n -\sigma}+G_7|h|^{2-\sigma}\Bigg]\\
		+\Bigg[G_8|h|^{2+\sigma n-2\sigma}+G_9|h|^{2+2\sigma n-2\sigma}+G_{10}|h|^{2-2\sigma}\Bigg]
	\end{eqnarray*}
	\begin{equation}\label{eq5.36}
		\hspace{110mm}+\Bigg[G_{11}|h|^{\gamma_1+\frac{2\sigma(\gamma_1-1)}{\gamma_1}}+G_{12}|h|^{\frac{2\sigma(\gamma_1-1)}{\gamma_1}+\sigma n}+G_{13}|h|^\frac{2\sigma(\gamma_1-1)}{\gamma_1}\Bigg]
	\end{equation}
	Where $G_1=\mathcal{H}_2+\mathcal{H}_3$ and $G_i=\mathcal{H}_{i-2}$ for $i=4,5,6,...,15$, are all finite terms.
	Since $|h|\leq 1$, then I have
	\begin{equation}\label{eq5.37}
		|h|^{2+2\sigma n-2\sigma},|h|^{2+\sigma n-2\sigma}\leq |h|^{2-2\sigma}
	\end{equation}
	\begin{equation}\label{eq5.38}
		|h|^{2+\sigma n-\sigma},|h|^{2+2\sigma n -\sigma},|h|^{2-\sigma},|h|^{2+\sigma n-2\sigma},|h|^{2+2\sigma n-2\sigma}\leq |h|^{2-2\sigma}
	\end{equation}
	and 
	\begin{equation}\label{eq5.39}
		|h|^{\gamma_1+\frac{2\sigma(\gamma_1-1)}{\gamma_1}}, |h|^{\frac{2\sigma(\gamma_1-1)}{\gamma_1}+\sigma n}\leq |h|^\frac{2\sigma(\gamma_1-1)}{\gamma_1}
	\end{equation}
	Now choosing $\sigma=\frac{1}{2}$ , then (\ref{eq5.36}) becomes
	\begin{equation}\label{5.40}
		\int_\mathfrak{B}|\tau_hV_{\gamma_1}(Du)|^2dx\leq c|h|+c|h|^\frac{\gamma_1-1}{\gamma_1}
	\end{equation}
	At this point, considering as in \cite{KM33}, we can start covering $\Omega'$ with the balls $\mathfrak{B}_1=\mathfrak{B}_1(x_1,|h|^\sigma),...,\mathfrak{B}_M=\mathfrak{B}_M(x_M,|h|^\sigma)$, where $M=M(h)\in \mathbb{N}$, such that corresponding inner cubes $\mathcal{Q}_1(\mathfrak{B}_1),...,\mathcal{Q}_1(\mathfrak{B}_M)$ are disjoints such that the sidelength of every such cube is comparable to $|h|^\sigma$ via a constant independent to $M$ and satisfy
	\begin{equation*}
		\Bigg|\Omega'\setminus\bigcup_{k=1}^{M}\mathcal{Q}_1(\mathfrak{B}_k)\Bigg|=0
	\end{equation*}
	Summing up (\ref{eq5.36}) with  inner cubes $\mathcal{Q}_1(\mathfrak{B}_1),...,\mathcal{Q}_1(\mathfrak{B}_M)$(since all constants will be independent of $h$, so patching up all the constants to some fixed constants),with $\sigma=\frac{1}{2}$ and with the help of (\ref{eq5.37},\ref{eq5.38},\ref{eq5.39}) we will have 
	\begin{equation}\label{eq197}
		\int_{\Omega'}|\tau_hV_{\gamma_1}(Du)|^2dx\leq c|h|+c|h|^\frac{\gamma_1-1}{\gamma_1}
	\end{equation}
	Now , integrating with respect to the measure $\frac{dh}{h^n}$ over the ball $B_\frac{r}{4}(0)$ , we obtain the following estimate 
	\begin{equation*}
		\int_{B_\frac{r}{4}(0)}\Bigg(\int_{\Omega'}\frac{|\tau_hV_{\gamma_1}(Du)|^2}{|h|^{2\theta}}dx\Bigg)^\frac{q}{2}\frac{dh}{|h|^n}\leq c\int_{B_\frac{r}{4}(0)}|h|^{(1-\theta)q}\frac{dh}{|h|^n}+c\int_{B_\frac{r}{4}(0)}|h|^{(\frac{\gamma_1-1}{\gamma_1}-\theta)q}\frac{dh}{|h|^n}
	\end{equation*}
	\begin{equation}\label{eq5.42}
		=T_1+T_2
	\end{equation}
	for $\theta<\frac{1}{\gamma_1}$.\vskip 1mm
	Now consider the integral $T_1$ and $T_2$, then by polar cordinates and assumption on $q$,
	\begin{equation*}
		T_1=c\int_{B_\frac{r}{4}(0)}|h|^{(1-\theta)q}\frac{dh}{|h|^n}=c\int_{0}^{r/4}\rho^{(1-\theta)q-1}\;d\rho<+\infty
	\end{equation*}
	\begin{equation*}
		T_2=c\int_{B_\frac{r}{4}(0)}|h|^{(\frac{\gamma_1-1}{\gamma_1}-\theta)q}\frac{dh}{|h|^n}=c\int_{0}^{r/4}\rho^{(\frac{\gamma_1-1}{\gamma_1}-\theta)q-1}\;d\rho<+\infty
	\end{equation*}
	where $T_2$ is finite because for $\gamma_1>2$ and $\theta<\frac{1}{\gamma_1}$, $(\frac{\gamma_1-1}{\gamma_1}-\theta)>\frac{\gamma_1-2}{\gamma_1}$, since by assumption $q>\frac{\gamma_1}{\gamma_1-2}$, so the integral is finite.\vskip 1mm
	Then by Lemma \ref{L3.9}, we can say from (\ref{eq5.42}) that $V_{\gamma_1}(Du)\in B^\theta_{2,q}(\Omega)$ locally.\hspace{10mm} $\square$\vskip 5mm
	{\large \bf Acknowledgement : } Author is grateful to  his teacher {\bf Dr. Ramesh Mondal} ( {\em Assistant Professor, Dept. of Mathematics, University of Kalyani}) for supporting throughout the work and giving some valuable inputs to improve the manuscipt. Author is thankful to {\bf CSIR-MHRD}, {\em Govt. of India} for helping the author with funds ({\em CSIR JRF} ) under file no. {\bf 09/0106(13571)/2022-EMR-I}.
	 \vskip 10mm
\textsc{Department of Mathematics\\University of Kalyani\\Kalyani, Nadia, West Bengal 741235\\India}\\[2mm] Email : {kardebraj98@gmail.com} , {debrajmath22@klyuniv.ac.in}
\end{document}